\documentclass[reqno,11pt]{amsart}
\usepackage{amsmath}
\usepackage[left=2.5cm,right=2.5cm,top=2.5cm,bottom=2.5cm]{geometry}
\usepackage{mathtools}
 \usepackage[abs]{overpic}

\usepackage[colorlinks=true, pdfstartview=FitV, linkcolor=blue, 
            citecolor=blue, urlcolor=blue]{hyperref}

\numberwithin{equation}{section}
\theoremstyle{plain}
\newtheorem{theorem}{Theorem}[section]
\newtheorem{proposition}[theorem]{Proposition}
\newtheorem{lemma}[theorem]{Lemma}

\newtheorem{definition}[theorem]{Definition}

\newcommand\R{\mathbb R}

\newcommand\N{\mathbb N}

\newcommand\wk{\rightharpoonup}
\newcommand\wstarto{\stackrel{\ast}{\wk}}
\newcommand\di{\mathrm{d}}

\newcommand\UUU{\color{black}}

\newcommand\EEE{\color{black}}
\newcommand\RRR{\color{black}}
\newcommand\MMM{\color{black}}

\title[Existence results for a morphoelastic growth model] {Existence results for a morphoelastic model}
\author[E. Davoli] {Elisa Davoli} 
\address[Elisa Davoli]{Institute of Analysis and Scientific Computing, TU Wien, 
Wiedner Hauptstrasse 8-10, A-1040 Vienna, Austria}
\email{elisa.davoli@tuwien.ac.at}
\author[K. Nik] {Katerina Nik}
\address[Katerina Nik]{Faculty of Mathematics, University of Vienna, 
Oskar-Morgenstern-Platz 1, A-1090 Vienna, Austria}
\email{katerina.nik@univie.ac.at}
\author[U. Stefanelli] {Ulisse Stefanelli} 
\address[Ulisse Stefanelli]{Faculty of Mathematics, University of Vienna, 
Oskar-Morgenstern-Platz 1, A-1090 Vienna, Austria,\,\&
Vienna Research Platform on Accelerating Photoreaction Discovery, University of Vienna, W\"ahringerstrasse 17, A-1090 Vienna, Austria,\,\&
Istituto di Matematica Applicata e Tecnologie Informatiche E. Magenes,
via Ferrata 1, I-27100 Pavia, Italy.
}
\email{ulisse.stefanelli@univie.ac.at}


\begin{document} 
\vskip .2truecm
\begin{abstract} We present some existence results for
  three-dimensional quasistatic 
  morphoelasticity. The state of the growing body is
  described by its deformation and the underlying growth tensor and is
  ruled by the interplay of hyperelastic energy minimization and
  growth dynamics. By introducing a regularization in the model, we prove that
  solutions can be obtained as limits of time-discrete solutions,
  built by means of an exponential-update scheme.
By further allowing
  the dependence of growth dynamics on an additional scalar field, to
  be interpreted as a nutrient or inhibitor, we formulate an optimal
  control problem and prove existence of optimal controls and
  states. Eventually, we tackle the existence of coupled morphoelastic
  and nutrient solutions, when the latter is allowed to diffuse and
  interact with the growing body. 
\end{abstract}
\maketitle

\section{Introduction}\label{sec1}
Morphoelasticity describes the growth of an elastic body and
finds its main application in the context of biological systems. Here,
growth is often a central aspect and is driven by a variety of
phenomena
acting at different scales. Below, we limit ourselves in summarizing some
key modeling issues, referring to the recent monograph \cite{goriely} for a thorough
introduction to the topic and additional material. 

The description of the
mechanical response of a growing body can be simplified by
restricting the attention to the macroscopic level of continua.
Assume to be given a nonempty, open, simply connected and bounded set
$\Omega\subset \R^3$ with smooth boundary, to be interpreted as the  {\it reference
configuration} of the body. At all times $t \in [0,T]$,  $T>0$, the {\it deformation}
of the body will be denoted by $y(t):\Omega \to \R^3$. 
Classical morphoelastic models postulate the {\it
  multiplicative decomposition} of the deformation gradient
$\nabla y(t)$  into an {\it elastic strain tensor} $F_{\rm el} \in
\R^{3\times 3}$,
related to stresses, and a {\it growth tensor} $G  \in
\R^{3\times 3}$, specifying the growth dynamics, namely, 
$$\nabla y(t)=F_{\rm el}(t)G(t).$$
In case $G(t)$ is compatible, namely, if $G(t) = \nabla y_{\rm
  gr}(t)$ for some given {\it growth deformation} $y_{\rm
  gr}(t):\Omega \to \R^3$, one can prove that $F_{\rm el}(t)$ is compatible
as well and the latter multiplicative decomposition
corresponds to the classical
chain rule applied to the composition $y(t ) = y_{\rm el}(t) \circ y_{\rm
  gr}(t)$. Here, $y_{\rm el}(t)$ can be
interpreted as the elastic deformation of the {\it evolved}
configuration $ y_{\rm gr}(t,\Omega)$. We however do not assume
compatibility here, for this would limit the applicability of the
theory, see \cite[Sec. 12.5]{goriely}. 

The {\it state} of the morphoelastic system is hence determined by the
pair $(y(t), G(t))$ for   $t\in [0,T]$.  
Its evolution in time is governed by the interplay
between the mechanical equilibrium and the growth process. As the time scales of
mechanical equilibration and of growth usually differ by 
orders of magnitude, inertial effects can be assumed to be negligible and one resorts
to the quasistatic approximation of the equilibrium system 
\begin{align}
  {\rm div}\,P(t)+ f(t)=0&\quad
                                                                 \text{in}
                                                                 \ \
                                                                 [0,T]\times
                                                                 \Omega,\label{eq:q1}\\
  y(t)= {\rm id}&\quad
                                                               \text{in}
                                                                 \ \
                                                                 [0,T]\times
                                                                 \Gamma_{\rm
                  D},\label{eq:q2}\\
  P(t) n = g(t) &\quad
                                                                 \text{in}
                                                                 \ \
                                                                 [0,T]\times
                                                                 \Gamma_{\rm
                  N} \label{eq:q3}.
\end{align}
The tensor $P(t)$ above is the {\it first Piola-Kirchhoff} stress. 
We assume the body to be {\it hyperelastic}, so that its elastic
state is determined by the elastic energy density $W=W(F_{\rm el}) =
W(\nabla y(t) G^{-1}(t))$. In particular,  $P(t)$ in \eqref{eq:q1}
is given by 
\begin{equation}\label{eq:Piola}
  P(t)= \UUU \det G(t)  \EEE\,{\rm D}W(\nabla y(t)G^{-1}(t))G^{-\top}(t).
  \end{equation}
The quasistatic equilibrium system   features a  time-dependent body force $f(t):
\Omega \to \R^3$ as well as a time-dependent surface traction $g(t):
\Gamma_{\rm N} \to \R^3$, to be imposed on the Neumann part
$\Gamma_{\rm N}$ of the boundary $\partial \Omega$. In addition, the body is clamped at
$\Gamma_{\rm D}\subset \partial \Omega$, where
$\Gamma_{\rm D}\cap \Gamma_{\rm N} = \emptyset$. 

The evolution of the growth tensor $G(t)$ is specified via the space-parametrized ODE in
rate form $$ G'(t) G^{-1}(t)=M(t)$$ where the prime stands for
partial time
differentiation.  The constitutive choice for the {\it growth
  rate}   $M$ reflects the 
combination of different effects driving the evolution and we refer the reader to \cite{Epstein,goriely} for a
discussion on the many possibilities. In all
generality, $M$ can depend on time $t\in [0,T]$, referential position
$x\in \Omega$, and actual position $y(t,x)$, modeling indeed nonhomogeneous
growth conditions in time and space. Growth may also be influenced by the
state of the system, namely, by $G(t)$ and by the deformation gradient
$\nabla y (t)$. In addition, the stress $P(t)$ is known to be
possibly driving
growth in some applications \cite{Jones}. 

 In the
following, we hence resort in focusing on a some reduced 
evolution model by
prescribing
\begin{align}
  \label{eq:q4}
   G'(t) G^{-1}(t) = M(G(t), \nabla y(t))&\quad \text{in} \ \
                                              [0,T]\times \Omega,\\
  G(0)=G^0 &\quad \text{in} \ \
                                              \Omega\label{eq:q5}
\end{align}
In the latter, all nonhomogeneities are neglected for the sake of
simplicity. Note however that these could be considered as well, at
the price of some additional notational intricacy. Notably, the
dependence on the stress $P(t)$ can be accounted for in
\eqref{eq:q4} by means of the dependence on the tensors $G(t)$
and $ \nabla y(t)$, by implicitly assuming
\eqref{eq:Piola}. Note in passing that the actual dependence of $M$ on   
stress or strain is still debated \cite{Ambrosi}.

 Research in morphoelasticity has
been up to now primarily devoted to clarifying the mechanical setting
and to deriving numerical simulations. In this respect, we refer the
reader to the recent
\cite{Chenchiah,Erlich, Haas} and \cite{Balbi,Fok,Kadapa}. To the best
of our knowledge, an existence
theory for solutions of the nonlinear morphoelastic evolution system
\eqref{eq:q1}-\eqref{eq:q5} is still
unavailable.

 In this paper, we move first steps in this direction,
by focusing on some nonlocally relaxed versions of the growth-dynamics
rule \eqref{eq:q4}. Indeed, the analysis of problem \eqref{eq:q1}-\eqref{eq:q5} requires formulating first a time-discrete version of the system, introducing suitable piecewise continuous and piecewise affine interpolants of the key quantities, and eventually passing
to the limit as the width of the time step tends to zero. In particular, a key point is the limit passage in the nonlinear rate $M$, which we will assume to be Lipschitz
continuous with respect to its variables. This in turn calls for some
time-compactness of  the interpolants of $\nabla y$, which however is not to be expected in
the quasistatic framework of \eqref{eq:q1}-\eqref{eq:q3}. In fact,
the best one can hope for from minimality is a uniform Sobolev bound,
see \eqref{eq:unif-y} later on. Consequently,
the analysis of \eqref{eq:q4} would soon grind to a halt. We hence propose to introduce a regularization of
$\nabla y$ in the dependence of $M$. This is achieved by 
replacing \eqref{eq:q4} by 
\begin{equation}
  \label{eq:q6}
  G'(t) G^{-1}(t) = M(G(t), (K \nabla y)(t))\quad \text{in} \ \
                                              [0,T]\times \Omega\\
\end{equation}
 where $ (K \nabla y)(t)$ is defined as a space and time
convolution as
\begin{equation}
  \label{eq:convolution}
  (K\nabla y)(t,x) = \int_0^t \int_{\R^3} \kappa (t-s)\phi(x-z)
  \nabla y(s,z)\, {\rm d} z \, {\rm d} t\quad \forall \, (t,x)\in [0,T]\times \Omega.
\end{equation}
Note that, here and in the following, $\nabla y$ is tacitly intended to be
trivially extended to zero in $\R^3\setminus \Omega$ whenever needed. 
This regularization serves the mathematical purpose of allowing a
satisfactory existence theory, at the price of a minor modification
to the model. In fact, the presence of the convolution term may be
justified from the modeling viewpoint as of  introducing a
nonlocal in space and time dependence of the rate $M$ on the
deformation gradient. In particular, the presence of the 
time-convolution 
kernel $\kappa$ induces a second time scale into the problem, which
may be interpreted as a time relaxation. Note that a similar effect
would have been achieved by considering viscoelastic dynamics
instead. 
 
Our first main result is an existence theory for a variational
formulation on the regularized
{\it morphoelastic evolution problem} \eqref{eq:q1}-\eqref{eq:q3},
\eqref{eq:q5}-\eqref{eq:q6}. As already mentioned, this is achieved via a
time-discretization argument. A
crucial observation here is that the sign of the determinant of $G(t)$
is preserved along the evolution, as an effect of the nonlinear structure of
\eqref{eq:q6}. This conservation is crucial, for it guarantees that
the growth process is nondegenerate and locally orientation preserving.
Correspondingly, we resort to a time-discrete scheme of
exponential type, reproducing this sign conservation at the
discrete level.

The structure of the morphoelastic evolution problem is
reminiscent of the quasistatic evolution problem in creep inelasticity, the difference being that in this latter setting $M$
is taken to be the variation with respect to $G$ of the total
energy functional. As such, in creep inelasticity $M$ is directly related to $W$. This
entails conservation of energy, at least
formally, which in turn provides the fundamental a priori
estimate. The present situation is different, for energy cannot
be expected to be conserved along the growing process, as we are not
including all energy exchanges in the description. In particular, no relation is imposed 
between the elastic-energy density  $W$ and the growth-rate function $M$, the specification of the latter being usually
just phenomenological. As a consequence, we have to obtain a priori estimates
otherwise. 

Note that most existence theories in
multidimensional inelasticity at finite strains hinge on the
presence of higher-order gradients in the internal variable $G$, here
to be interpreted as inelastic strain
\cite{cplas_part2,mainik-mielke2,Mielke-Mueller,mielke-roubicek16}
(see \cite{compos2,RoegerSchweizer17}, however, where no gradient is
involved). 
We
avoid such higher-order terms here, still allowing some nonlocal
effect in space
via the convolution term $K\nabla y$.

Our existence result can be compared with the one in
\cite{MRS_vis}. There, a similar model to \eqref{eq:q1}-\eqref{eq:q5}
is introduced in the frame of rate-dependent
viscoplasticity and proved to admit solutions via a
time-discretization and passage to the limit procedure. The existence theory in
\cite{MRS_vis} is however quite different from ours. At first, the
analysis in \cite{MRS_vis} hinges upon assuming a variational origin
of the flow rule, which is not avaliable here. Secondly, a gradient term
in $G$ is considered, whereas our model is local in $G$. Thirdly, the
solution notion in \cite{MRS_vis} is variational, making the
regularization of the occurrence $\nabla y$ in the flow rule unnecessary. Eventually,
the time-discretization scheme in \cite{MRS_vis} is the classical
variational one, while we consider an {\it exponential} variant instead,
see \eqref{eq:discrete-growth2}.

As a second existence result, we consider an additional dependence of $M$ 
 from an external field $\mu(t)$, which we assume to
be scalar for definiteness. Namely, we replace \eqref{eq:q6} by 
\begin{equation}
  \label{eq:q7}
  G'(t) G^{-1}(t) = M(G(t), (K\nabla y)(t),\mu(t))\quad \text{in} \ \
                                              [0,T]\times \Omega.\\
\end{equation}
The field $\mu(t) : \Omega \to \R$ can be interpreted as the
concentration of a nutrient (or an inhibitor), influencing the growth
rate. In Section \ref{sec:control} we analyze an optimal control
problem, where $\mu(t)$ acts as a control and drives the trajectory
$t \mapsto (y(t),G(t))$ to minimize a target functional, possibly of the form 
$$J(y,G) = \int_\Omega i(G(T))\, {\rm d} x + \int_0^T\int_\Omega
j(y,G)\, {\rm d} x \, {\rm d} t.$$
Here, the functions $i$ and $j$ are suitably lower semicontinuous and favor
specific deformations and growth tensors. In particular, the choice
$i(G) = \det G$, $j=0$ represents volume minimization whereas
$j(y,G) =  |y-y_{\rm target}|^2$, $i=0$   corresponds to the possible
attainment of a given target deformation $y_{\rm target}$.

Eventually, we couple the evolution of the state $(y(t),G(t))$ with
that of the scalar field $\mu(t)$ by additionally specifying its
evolution as 
\begin{align}
 \mu'(t) - \nu \Delta \mu(t) = h(t) - H( (\kappa \ast y)  (t))& \quad \text{in} \ \
                                              [0,T]\times \Omega,  \label{eq:mu1}\\
 \mu(t) = \mu_{\rm D}(t) &  \quad \text{in} \ \ 
                                              [0,T]\times  \partial \Omega,  \label{eq:mu2}\\
 \mu(0) = \mu^0  &  \quad \text{in} \ \ 
                                               \Omega  \label{eq:mu3}. 
\end{align}
Here, $\nu>0$ and $h(t):\Omega \to \R$ plays the role of a given
source. The
term $H(\kappa \ast y )$ instead models the consumption of $\mu(t)$ during the growth
process, which is indeed assumed to depend on the actual position of
the body, again mollified by a time-convolution compactifying
term. 
 In particular, the triplet $(y(t),G(t),\mu(t))$ describes a
system, where growth is influenced by the field $\mu$
which diffuses and is consumed during growth.  The coupling of the
quasistatic equilibrium \eqref{eq:q1}-\eqref{eq:q3}, the growth
dynamic \eqref{eq:q5},
\eqref{eq:q7}, and the nutrient dynamic \eqref{eq:mu1}-\eqref{eq:mu3}
gives rise to a {\it nutrient-morphoelastic evolution
  problem}, which is variationally reformulated and proved to admit
solutions in   Section \ref{sec:chemotaxis} below. Prototypical phenomena encoded by the above system are those in which diffusion happens on a 
much slower time scale with respect to that of mechanical equilibration.  We hence keep track here of viscous effects in the nutrient dynamics.
For completeness, we mention that an alternative modeling choice would be to replace \eqref{eq:mu1} by a quasistatic evolution of the nutrient as well. 
This latter scenario could still be included in our analysis at the mathematical price of introducing a further nonlocality in the dependence of $M$ on $\mu$.

 The paper is organized as follows. In Section \ref{sec:setting} we
 introduce the precise mathematical framework and state our main
 results. Section \ref{sec:analysis} is devoted to the proof of
 Theorem \ref{thm:main}, in which concentration of nutrients is
 neglected. This latter dependence is accounted for in the control
 problem formulated in Theorem \ref{thm:control} whose proof is the
 subject of Section \ref{sec:control}. Eventually, the full
 nutrient-morphoelastic evolution problem is analyzed in Section \ref{sec:chemotaxis}.

\section{Setting and main results}
\label{sec:setting}

We devote this section to making assumptions precise and stating our
existence results. As anticipated in the Introduction, we
let the reference configuration of the body $\Omega\subset \R^3$ be
nonempty, open,  simply connected, bounded, and smooth and 
$\Gamma_{\rm N},\,\Gamma_{\rm D}\subset \partial \Omega$, with
$\Gamma_{\rm D}$ and $\Gamma_{\rm N}$ open in the
topology of $\partial \Omega$ and disjoint, $\Gamma_{\rm D}\not =
\emptyset$, and
$ \overline \Gamma_{\rm D} \cup \overline \Gamma_{\rm N} = \partial
\Omega$.

Throughout the paper, $GL_+(3)$ and $SO(3)$ denote the general
linear group and the set of proper rotations, i.e., 
\begin{align*}
 GL_+(3)   =\{ A \in \R^{3\times 3}:\,{\rm det}\, A>0\},
 \quad 
SO(3) =\{ A \in \R^{3\times 3}:\,{\rm det}\, A=1,\, A^\top A={\rm Id}\}, 
\end{align*}
where $\top$ denotes transposition and $\rm Id$ is the identity
$2$-tensor. Given the $2$-tensors $A,B 
\in \R^{3\times 3}$ and the $3$-tensors $C,D \in  \R^{3\times 3 \times 3}$ we classically define 
$A : B, C:D \in \R$ and $C:B, B:C \in \R^{3}$  as (summation convention) $A :B  := A_{ij}B_{ij}$, $C :D := C_{ijk}D_{ijk}$, $(C:B)_i := C_{ijk}B_{jk}$, and $(B:C)_i := B_{jk}C_{jki}$, respectively. The space of $2$-tensors $\R^{3\times 3 }$ is endowed with the natural scalar product $A:B:= {\rm tr\, }(A^\top B)$, where 
${\rm tr\,}(A):= A_{ii}$  and corresponding norm $|A|^2:= A:A$. We note that this norm is submultiplicative, i.e., $|AB| \leq |A| |B|$. Similarly, we define the norm $|C|^2:= C:C$, the partial transposition 
$(C^t)_{ijk}:= C_{jik}$, and the products $(CB)_{ijk}:= C_{ijl}B_{lk}$ and $(BC)_{ijk}:= B_{il}C_{ljk}$ so that 
$|CB|, |BC| \leq |B||C|$. 

 Furthermore, for $A,B: \Omega \to \R^{3 \times 3}$ 
 differentiable, the gradient $\nabla A$ is a 
$3$-tensor reading $(\nabla A)_{ijk}:=A_{ij,k}$ and it holds that $\nabla (AB)= (B^\top \nabla A^\top)^t + A 
\nabla B$. 

We will denote by ${\rm id}$ the identity map ${\rm id}(x)=x$ for all
$x\in \R^3$. Given a time-dependent map $\psi(t)$, we indicate
by ${\psi}'(t)$ its (possibly partial) derivative with respect to time.
 
 In all of the following, we assume to be given 
$$p>3$$
and denote by $q =p/ (p-1)<3/2$ the corresponding conjugate exponent.
Let us define the class of admissible deformations $\mathcal{Y}$
and admissible growth-tensor fields $\mathcal{G}_\infty$ as  
\begin{align}
	\label{eq:def-class-y}
\mathcal{Y}& := \{y\in W^{1,p}(\Omega; \R^3):  \,{\rm
             det}\,\nabla y>0\ \text{ a.e. in }\Omega ,
             \ y={\rm id}\ \text{ on }\Gamma_{\rm D}\},\\
\label{eq:def-class-G}
\mathcal{G}_\infty&:=\{G\in W^{1,\infty}(\Omega;\R^{3\times 3}):\,{\rm
             det}\,G>0\ \text{ a.e. in }\Omega\}.
\end{align}
The prescription on the a.e. positivity of $\det G$ in $\mathcal{G}_\infty$
is intended to guarantee that $G$ is not degenerate and is orientation
preserving. 
 
For a given growth tensor field $G(t) \in \mathcal{G}_\infty$, the
variational formulation of the quasistatic equilibrium system
\eqref{eq:q1}-\eqref{eq:q3} corresponds to the minimization on
$\mathcal{Y}$ of the {\it total elastic energy}
\begin{equation}
\label{eq:en}\mathcal{E}(y,G(t)) :=  \int_{\Omega} W(\nabla
y(x)G^{-1}(t,x))\,\UUU \det G(t,x)\EEE\,\di x-  \langle \ell(t), y \rangle,
\end{equation}
where we have indicated by $\ell(t) \in (W^{1,p}(\Omega;\R^3))'$
(dual) the
{\it generalized load} 
$$\langle \ell(t) , y \rangle := \int_\Omega f(t)\cdot y \, {\rm d} x
+ \int_{\Gamma_{\rm N}} g(t) \cdot y \, {\rm d} \mathcal{H}^2. $$
Here, $\langle \cdot, \cdot \rangle$ denotes the duality pairing between $(W^{1,p}(\Omega;\R^3)) '$ and 
$W^{1,p}(\Omega;\R^3)$ and $\mathcal{H}^2$ is the two-dimensional
Hausdorff measure. \UUU The explicit occurrence of $ \det G(t,x)
$ in the total elastic energy is a consequence of the fact that the
integration is taken with respect to the pre-growth reference
configuration $\Omega$ \cite{Kadapa,Moulton20}. \EEE

Above, 
$W:\R^{3\times 3}\to [0,\infty]$ denotes the {\it elastic energy
density}. We assume that $W \in C^1(GL_+(3))$, $W\equiv \infty$ on
$\R^{3\times 3}\setminus GL_+(3) $ and that it satisfies the following standard hypotheses:
\begin{itemize}
\item[\bf{(H1)}] (polyconvexity)  $\exists \widehat{W} : \R^{3\times 3}
  \times \R^{3\times 3} \times \R \to [0,\infty]$ convex and such that
$$
	W(A)= \widehat{W}(A, {\rm cof} \, A, {\rm det}\, A)  \quad \forall \, A
        \in GL_+(3). 
$$
\item[\bf{(H2)}] (coercivity and control)    $\exists c_1,\,
  c_2>0$ such that
  \begin{align*}
  W(A) \geq c_1|A|^p - \frac{1}{c_1}    \quad \text{and} \quad 
   |A^\top \partial_A W(A)| \leq c_2(W(A) +1)  \quad \forall \, A \in GL_+(3).
  \end{align*}
\end{itemize}
The second assumption in {\bf(H2)} prescribes the controllability of
the {\it Mandel tensor} $A^\top \partial_A W(A)$ via the
energy \cite{Bal84b,Bal02} and turns out to be particularly relevant in connection
with finite-strain elastoplasticity \cite{FM06,mainik-mielke2,ms5}. 
Note that assumptions {\bf (H1)}-{\bf (H2)} are compatible with
frame-indifference, namely, 
$$
W(RA)=W(A)\quad \forall  \, R\in SO(3),\   \forall \, A
\in GL_+(3). 
$$
We will check below that ${\rm det}\, G >0$ for all times. As $W$ is
unbounded out of $GL_+(3)$ only, we hence have that ${\rm det}\, \nabla y
>0$ for a.e. times, as soon as the energy is finite. 

Concerning body forces and traction we assume 
\begin{itemize}
\item[\bf{(H3)}] 
 $ f \in W^{1,1}(0,T; L^1(\Omega; \R^3))$ and $ g \in W^{1,1}(0,T; L^1 (\Gamma_N; \R^3))$.
\end{itemize}
Note in particular, that {\bf (H3)} entails 
$$
\ell \in W^{1,1}(0,T; (W^{1,p}(\Omega; \R^3)) ').
$$
 Let us recall that the evolution of $G(t)$ is driven by 
relation  \eqref{eq:q6}, featuring the convolution term $K \nabla
y$.  To this effect, we specify 
\begin{itemize}
\item[\bf{(H4)}] $\kappa\in W^{1,1}(0,T)$ and $\phi\in W^{1,q}(\R^3)$
\end{itemize}
and define the operator $K$ as 
$$
(K \psi)(t,x):=(\kappa \ast (\phi \star \psi))(t,x)
\quad \forall \, \psi\in L^1((0,T)\times \R^3; \R^{3\times 3}), 
$$
where $\ast$ and $\star$ denote the standard convolution products on $(0,t)$ and $\R^3$, respectively. Namely,
$$
(\kappa \ast \psi)(t,\cdot):= \int_0^t \kappa (t-s) \psi(s,\cdot)\,\di s \quad \text{ for } t\in (0,T)
$$
and 
$$
(\phi \star \psi)(\cdot,x):= \int_{\R^3} \phi(x-z)\psi(\cdot,z)\,\di z \quad \text{ for } z \in \R^3. 
$$
By applying $K$ to (components of)
functions defined on $\Omega$ only, we actually consider the
corresponding trivial extensions to zero to the whole $\R^3$, without
introducing new notation. As regards the initial values, we assume 

\begin{itemize}
	\item[\bf{(H5)}]  $G^0\in \mathcal{G}_\infty$, \UUU  $\det
          G^0 \geq \delta  $ a.e. for some $\delta>0$, \EEE and
          $y^0 \in \text{arg\,min}_{\mathcal{Y}}
          \mathcal{E}(\cdot,G^0)$. 
\end{itemize}
 Note that we will prove  in Lemma \ref{lemma:exist-y} below that
such a minimizer $y^0$ exists for all $G^0\in \mathcal{G}_\infty$. 

 We will work under the following regularity of the growth-rate
 function 
\begin{itemize}
	\item[\bf{(H6)}]  $M\in W^{1,\infty}( \R^{3\times 3}\times 
	\R^{3\times 3}; \R^{3\times 3})$. 
      \end{itemize}

From assumptions {\bf{(H5)}} and {\bf{(H6)}}, for all $t \mapsto
G(t)$ solving \eqref{eq:q6} it follows
that $\det G(t) >0$ a.e. in $\Omega$, $\forall \, t \in [0,T]$. 
Indeed, by the Jacobi formula and equation \eqref{eq:q6} we have that 
\begin{equation*}
	\frac{\di}{\di t} \det G(t) = \det G(t) \, {\rm tr \, } M\big( G(t), (K{\nabla y}) (t)\big).
\end{equation*}
Solving this ODE gives 
\begin{equation}
 \det G(t) = \det G^0 \, \exp  \left(\int_0^t  {\rm tr \, } M\big(
   G(s), (K {\nabla y}) (s)\big) \, {\rm d} s \right). \label{eq:nondeg}
\end{equation}
Using {\bf{(H5)}} and estimating the integrand above yields 
\begin{equation*}
	\det G(t) \geq  \det G^0 \, \exp\big(-3T \, \Vert M\Vert_{L^{\infty}}  \big) >0 \quad \text{ a.e. in } \Omega, \, 
	\forall  \, t \in [0,T], 
\end{equation*}
where, here and in the rest of the paper, we use the short-hand
$\Vert \cdot \Vert_{L^{\infty}}$ to identify any $L^\infty$ norm, in
this case $
\Vert \cdot \Vert_{L^{\infty}(\R^{3\times 3}\times \R^{3\times 3};
  \R^{3\times 3})}$. \UUU This lower bound on $\det G(t)$ will turn
out crucial in combination with the coercivity in {\bf{(H2)}} in order
to prove the coercivity of the total elastic energy $\mathcal E$. \EEE 

\begin{definition}[Morphoelastic solution] We say that $(y,G) :[0,T] 
  \to \mathcal{Y} \times \mathcal{G}_\infty$ is a \emph{morphoelastic solution} if
\begin{align}
\label{eq:problemP}  &y(t)\in {\rm arg\,min}_{y\in \mathcal{Y}}\,
  \mathcal{E}(y,G(t))\quad\text{for a.e.} \ t\in (0,T), \\
& G'(t)G^{-1}(t) = M(G(t), (K \nabla y)(t))
  \quad\text{a.e. in} \ \Omega, \ \text{for a.e.} \ t\in
  (0,T), \label{eq:problemG}\\
&(y(0),G(0)) = (y^0,G^0)  \quad\text{a.e. in} \ \Omega.\label{eq:initial}
\end{align}
\end{definition}

Our basic existence result is the following.
\begin{theorem}[Morphoelastic existence]
\label{thm:main} Under assumptions {\bf (H1)}-{\bf (H6)}  there
exists a morphoelastic solution 
~$(y,G)\in L^\infty(0,T;W^{1,p}(\Omega; \R^3))  \times L^\infty(0,T;W^{1,\infty}(\Omega; \R^{3\times 3})) 
\cap W^{1,\infty}(0,T;L^\infty(\Omega; \R^{3\times 3}))$.
\end{theorem}

 The proof of Theorem \ref{thm:main} is in Section
 \ref{sec:analysis} below.

 Let us now turn to the case where the growth dynamics is
 influenced by the given nutrient field $\mu: [0,T] \times \Omega  \to \R$. To this aim, the
 growth-rate function $M$ has to be modified by including an
 additional dependence on the nutrient field
 $\mu$. Also in this extended case, we assume $M$ to be Lipschitz
 continuous, namely,  we modify {\bf (H6)} as 
\begin{itemize}
	\item[\bf{(H7)}] $M\in W^{1,\infty}( \R^{3\times 3}\times 
	\R^{3\times 3} \times \R;  \R^{3\times 3})$.
\end{itemize}
Correspondingly, we specify the class of admissible growth tensor
fields as
$$
	\mathcal{G}_p:=\{G\in W^{1,p}(\Omega;\R^{3\times 3}):\,{\rm
		det}\,G>0\quad\text{a.e. in }\Omega\}.
$$
One can define the following.

\begin{definition}[Nutrient-driven morphoelastic solution] 
Assume to be given  $\mu \in L^p(0,T;W^{1,p}(\Omega))$. We say that $(y,G) :[0,T] 
  \to \mathcal{Y} \times \mathcal{G}_p$ is  a
  \emph{nutrient-driven morphoelastic solution given $\mu$}  if
\begin{align}
\label{eq:problemPn}  &y(t)\in {\rm arg\,min}_{y\in \mathcal{Y}}\,
  \mathcal{E}(y,G(t))\quad\text{for a.e.} \ t\in (0,T), \\
& G'(t)G^{-1}(t) = M(G(t), (K {\nabla y})(t),\mu(t))
  \quad\text{a.e. in} \ \Omega, \ \text{for a.e.} \ t\in
  (0,T), \label{eq:problemGn}\\
&(y(0),G(0)) = (y^0,G^0)  \quad\text{a.e. in} \ \Omega.\label{eq:initialn}
\end{align}\label{nutrient-driven} 
\end{definition}

In Section~\ref{sec:control} we check that the existence result of
Theorem \ref{thm:main} can be readily extended to include the
nutrient-driven case. In particular, one can define a possibly set-valued
solution operator
\begin{align*}
&S: L^p(0,T;W^{1,p}(\Omega)) \\
&\quad\to L^\infty(0,T;W^{1,p}(\Omega; \R^3))  \times 
 W^{1,\infty}(0,T;L^\infty(\Omega; \R^{3\times 3})) \cap W^{1,p}(0,T;W^{1,p}(\Omega; \R^{3\times 3})) 
 \end{align*}
defining the set $S(\mu)$ of all nutrient-driven morphoelastic
solutions $(y,G)$ given $\mu$,
according to Definition \ref{nutrient-driven}. One can hence use the
solution operator $S$ to specify the optimal control problem
\begin{equation}
  \label{eq:oc}
  \min_{\mu \in \mathcal{A}} \big\{  J(y,G,\mu) \ : \ (y,G) \in S(\mu) \big\}.
\end{equation}
 Here, $\mathcal{A}\subset L^p(0,T;W^{1,p}(\Omega))$ is the set of
 admissible controls $\mu$. We assume that
 \begin{itemize}
	\item[\bf{(H8)}]  $\mathcal{A}$ is bounded in $
          L^p(0,T;W^{1,p}(\Omega))$ and compact in $L^1((0,T) \times\Omega)$. 
	\item[\bf{(H9)}]  $J: L^\infty(0,T;W^{1,p}(\Omega; \R^3))
          \times C([0,T];L^\infty(\Omega; \R^{3\times
            3}))\times L^p(0,T;W^{1,p}(\Omega)) \to  [0,\infty] $
          is lower semicontinuous with respect to the corresponding
          weak* topology.
        \end{itemize}
As already mentioned in the Introduction, this assumption on $J$
allows flexibility with respect to the possible choices
for $J$. These include, in particular, 
$$
J(y,G,\mu) = \beta_1\int_\Omega \det G(T)\, {\rm d} x +  \beta_2 \int_0^T\int_\Omega 
|y - y_{\rm target}|^p\, {\rm d} x \, {\rm d} t +  \beta_3 \int_0^T\int_\Omega |\mu|^p \, {\rm
  d} x \, {\rm d} t,
$$
which, together with incompressibility  (i.e., ${\rm det}\,
\nabla y =1$), would correspond to a
weighted combination ($\beta_1, \, \beta_2 ,\, \beta_3 \geq 0$)  of final
volume minimization, attainment of a given {\it target} deformation $y_{\rm target} \in
L^\infty( 0,T,W^{1,p}(\Omega;\R^3))$, and minimization of the amount of
provided nutrient (in connection with tumor
growth, one may think here of a chemotherapy drug).

Solutions to the optimal control problem \eqref{eq:oc} are optimal
controls $\mu^*$ and corresponding optimal pairs $(y^*,G^*)\in
S(\mu^*)$. Our next
result guarantees that these exist.

\begin{theorem}[Existence of optimal controls] \label{thm:control}
  Under assumptions {\bf(H1)}-{\bf(H5)}, {\bf(H7)}-{\bf(H9)} the
  solution operator $S$ is well-defined and the optimal control
  problem \eqref{eq:oc}
admits a solution $(y^*,G^*,\mu^*)$.
\end{theorem}

Theorem \ref{thm:control} is proved in 
Section \ref{sec:control}.

A further extension of the model corresponds to considering the driving
nutrient field to be unknown and to evolve together with the
mechanical variables, as effect of system \eqref{eq:mu1}-\eqref{eq:mu3}. To this aim, we introduce the 
class of admissible nutrient concentrations as
\begin{equation*}
\mathcal{M}:=W^{2,p}(\Omega)
\end{equation*}
and qualify boundary and initial data as
\begin{itemize}
	\item[\bf{(H10)}]   $\mu_{\rm D} \in  W^{1,\infty}(0,T;L^p(\Omega))\cap L^{\infty}(0,T;
          W^{2,p}(\Omega))$  and $\mu^0 \in \mathcal{M}$. 
\end{itemize}
The nutrient source and consumption terms in the right-hand side of equation \eqref{eq:mu1} are
provided via  
\begin{itemize}
	\item[\bf{(H11)}]  $h\in L^{\infty}( 0,T; L^p (\Omega))$ and
          $H  \in W^{1,\infty}(\R^3)$. 
        \end{itemize}
   Note in particular that we have $H (\kappa \ast y) \in L^{\infty}( 0,T; L^p
        (\Omega))$ whenever  $y\in L^{1 }( 0,T; L^p (\Omega; \R^3))$.
We are now ready to define our concept of solution of the fully
coupled system. 

\begin{definition}[Nutrient-morphoelastic solution] 
	  \label{nutrient-morphoelastic}
We say that $(y,G,\mu) :[0,T] 
\to \mathcal{Y} \times \mathcal{G}_p\times \mathcal{M}$ is a
\emph{ 
nutrient-morphoelastic solution}  if
\begin{align}
	\label{eq:problemcP}  &y(t)\in {\rm arg\,min}_{y\in \mathcal{Y}}\,
	\mathcal{E}(y,G(t))\quad\text{ for a.e. } \ t\in (0,T), \\
	& G'(t)G^{-1}(t) = M(G(t), (K \nabla y)(t), \mu(t))
	\quad\text{ a.e. in } \ \Omega, \ \text{ for a.e. } \ t\in
	(0,T), \label{eq:problemcG}\\
	&  \mu'(t) - \nu \Delta \mu(t) = h(t) - H((\kappa \ast y)(t)) \quad \text{ a.e. in } \ \Omega, \ \text{ for a.e. } \ t\in
	(0,T), \label{eq:problemcmu}\\
	&\mu(t) = \mu_{\rm D}(t)  \quad\text{ a.e. on} \ \partial \Omega,  \ \text{ for a.e. } \ t\in
	(0,T), \label{eq:boundaryc-1} \\
	& (y(0),G(0),\mu(0)) = (y^0,G^0,\mu^0) \quad\text{a.e. in} \ \Omega \label{eq:initialc-1}. 
\end{align}
\end{definition}

We are eventually in the position of presenting an existence result
for 
nutrient-morphoelastic solutions. 

\begin{theorem}[Nutrient-morphoelastic existence] \label{thm:chemotaxis}
Under assumptions {\bf (H1)}-{\bf (H5)}, {\bf (H7)},  and {\bf
  (H10)}-{\bf (H11)} there exists a  nutrient-morphoelastic solution $(y,G, \mu)\in L^\infty(0,T;W^{1,p}(\Omega; \R^3))  \times W^{1,\infty}(0,T;L^{\infty}(\Omega; \R^{3\times 3})) \cap W^{1,p}(0,T;W^{1,p}(\Omega; \R^{3\times 3})) \times 
W^{1,\infty}(0,T; L^p(\Omega)) \cap L^\infty(0,T; W^{2,p}(\Omega))$. 
\end{theorem}

Theorem \ref{thm:chemotaxis} is proved in Section \ref{sec:chemotaxis}.

\section{Proof of Theorem 1:  Morphoelastic existence}
\label{sec:analysis}

This section is devoted to prove the existence of morphoelastic
solutions, namely, trajectories $t \in
[0,T] \mapsto (y(t),G(t)) \in \mathcal{Y}\times \mathcal{G}_\infty$
fulfilling \eqref{eq:problemP}-\eqref{eq:initial}. 
We argue by time-discretization: we obtain time-discrete
solutions, prove a-priori estimates for the piecewise affine and
backward piecewise constant time-discrete interpolants, and eventually
pass to the limit as the  time step  converges to zero. For
convenience of the reader, each of the above steps is associated to a
corresponding subsection.

 In order to shorten notation, from here on we use the symbols ~$\Vert \cdot
\Vert_{L^{\infty}}$ and $\Vert \cdot \Vert_{W^{1,\infty}}$ to indicate generic
$L^\infty$ and $W^{1,\infty}$ norms, without explicitly specifying
dependencies. 

\subsection{Time discretization}
\label{subs:time-discrete}
We consider a uniform partition $\{ 0=t_0<t_1<\dotsc<t_{N-1}<t_N=T \}$, $N\in \N$, $t_i =i \tau$, $\tau >0$, 
of the time interval $[0,T]$. Given any vector $\{ z_i \}_{i=0}^N$, we will denote by $\hat{z}_{\tau}$ and 
$\bar{z}_{\tau}$ the corresponding piecewise affine and backward piecewise constant interpolants associated to the partition. Namely, 
\begin{align}
&\label{eq:pw-aff}\hat{z}_{\tau}(0):= z_0, \quad \hat{z}_{\tau}(t):= \alpha_i (t) z_i + (1-\alpha_i (t)) z_{i-1}, \\[0.1cm]
&\label{eq:pw-con}\bar{z}_{\tau}(0):= z_0, \quad \bar{z}_{\tau}(t):=z_i \qquad \text{ for } t \in ( (i-1)\tau, i\tau], \; i=1,\dotsc,N,
\end{align}
where $\alpha_i(t):= (t-(i-1) \tau)/ \tau$ for $ t \in (
(i-1)\tau, i\tau]$,  $i=1,\dotsc,N$.
Setting $\ell_i := \ell(t_i)$ for $i=0,\dots,N$, we obtain a discrete solution $\{ (y_i,G_i) \}_{i=1}^N \in \mathcal{Y}^N \times 
\mathcal{G}_\infty^N$ by recursively solving 
\begin{align}
&\label{eq:discrete-growth1} y_i\in \text{argmin}_{y \in
                \mathcal{Y}}\,\left\{\int_{\Omega} W(\nabla y \,
                G_{i}^{-1})\,\UUU \det G_{i} \EEE \,\di x- \langle \ell_i,y\rangle
\right\}, \\
&\label{eq:discrete-growth2} G_i = \exp \big(\tau \,  M(G_{i-1}, (K_\tau {\nabla y})_{i-1}) \big) \, G_{i-1} \quad \text{ a.e. in } \Omega, 
\end{align}
for $i=1,\dotsc,N$, starting from the initial data $(y_0,G_0)= (y^0,G^0) \in \mathcal{Y} 
\times \mathcal{G}_\infty$ \UUU with $\det G^0 \geq \delta>0$ a.e. \EEE As already mentioned, this {\it
  exponential-update} scheme is designed to reproduce at the discrete
level the nonlinear geometry of the differential system
\eqref{eq:problemG}, in particular the nondegeneracy property
\eqref{eq:nondeg}.

In \eqref{eq:discrete-growth2}, the operator $(K_\tau {\nabla y})_{i-1}$ is given by 
\begin{equation}
\label{eq:discrete-conv}
	(K _\tau {\nabla y})_{i-1}(x):= (\kappa \ast_\tau (\phi \star \nabla y) )_{i-1} (x) \quad\text{for a.e. } x \in \Omega,
\end{equation}
where the time discrete convolution $\ast_\tau$ is defined as 
$$
(\kappa  \ast_\tau \nabla y)_{i-1}:=  \sum_{j=0}^{i-1} \tau \,  \kappa_j \,{\nabla y}_{i-1-j}  \quad\text{ for } i=1,\dotsc,N
$$
with $\kappa_i:=\kappa(t_i)$, $ i=0,\dotsc,N$, \RRR see \cite{vol}. In order to prove
the existence of discrete solutions, we start by considering the
 equilibrium problem in the following lemma.

\begin{lemma}[Equilibrium problem]
  \label{lemma:elastic} Under assumptions \MMM{\bf (H1)}\EEE--{\bf (H2)}
  for every $G \in L^\infty(\Omega;\R^{3\times 3})$  with $\det G \geq \eta>0$
  a.e. for some $\eta>0$\MMM, and for \EEE every $\ell \in (W^{1,p}(\Omega;\R^{3}))'$ there exists $y\in \mathcal{Y}$ solving
  \begin{equation}
    \label{eq:elastic}
    \int_{\Omega} W(\nabla  y \,
                G ^{-1})\, \det G  \,\di x- \langle \ell , 
                y\rangle \leq 
    \int_{\Omega} W(\nabla \hat y \,
                G ^{-1})\, \det G  \,\di x- \langle \ell , \hat
                y\rangle \quad \forall \hat y \in \mathcal{Y}.
              \end{equation}
            \end{lemma}

            \begin{proof}
 Since $W \in C^1(GL_+(3))$, for all $\eta, \,g>0$  we have that
 \begin{equation}
   \MMM \lambda\RRR(g,\eta):= \max \{W(A^{-1}) \ : \ A \in GL_+(3), \ |A|\leq g, \
 \det A \geq \eta\}<\infty. \label{k}
 \end{equation}
Recalling that ${\rm id}\in \mathcal{Y}$ we get
 \begin{align}
   &\inf_{y \in \mathcal{Y}}\,\left\{\int_{\Omega} W(\nabla y \,
   G^{-1})\, \det G  \, \di x- 
 \langle \ell,y\rangle\right\}\leq \int_{\Omega} W(G^{-1})\,  \det G
\,\di x- \langle \ell,{\rm id}\rangle\nonumber\\
&\qquad\leq 6
|\Omega|\MMM \lambda\RRR(\|G\|_{L^{\infty}},\eta) \| G\|_{L^\infty}^3+ \| \ell\|_{(W^{1,p}(\Omega;\R^3))'}\|
{\rm id}\|_{W^{1,p}(\Omega;\R^3)} <\infty.\label{eq:stimaen}
\end{align}
Owing to the coercivity from {\bf (H2)}, every minimizing sequence
 $\{y_k\}\subset \mathcal{Y}$ fulfills
 \begin{align*}
 & \eta  \|\nabla y_k\|^p_{L^p(\Omega; \R^{3\times 3})} \leq  \eta  \|\nabla y_k \, G^{-1}\|^p_{L^p(\Omega; \R^{3\times 3})}
 \|G\|^p_{L^{\infty}} \leq  \frac{\| G\|_{L^{\infty}}^p}{c_1}\int_{\Omega}W(\nabla y_k \,
                                                  G^{-1})\,  \det G  \,\di{x}
                                                  +\frac{  \eta |\Omega|\| G\|_{L^{\infty}}^p}{c_1^2} \\
   &\quad \leq \frac{\| G\|_{L^{\infty}}^p}{c_1}\left(\int_{\Omega} W(G^{-1})\,  \det G
\,\di x+ \langle \ell,y_k -{\rm id}\rangle \right)   
     +\frac{ \eta |\Omega|\| G\|_{L^{\infty}}^p}{c_1^2}\\
   & \quad \leq  \frac{\| G\|_{L^{\infty}}^p}{c_1}\left(6  |\Omega| \MMM \lambda\RRR(\| G\|_{L^{\infty}},\eta) \| G\|_{L^{\infty}}^3 +\| \ell\|_{(W^{1,p}(\Omega;\R^3))'}\|y_k-{\rm id}\|_{W^{1,p}(\Omega;\R^3)}\right)+\frac{  \eta |\Omega|\|
     G\|_{L^{\infty}}^p}{c_1^2} \\
&\quad \leq C(1+\|\nabla y_k\|_{L^p(\Omega; \R^{3\times 3})})
 \end{align*}
for all $k \in \N$, where in the second-to-last inequality we have
used \eqref{eq:stimaen}, and where the last inequality
follows by {\bf (H3)} and the Poincar\'e inequality. The constant
$C>0$ depends on $\eta$, $c_1$, $|\Omega|$, $\| G\|_{L^{\infty}}$, $\MMM \lambda\RRR(\| G\|_{L^{\infty}},\eta)$, and $\| \ell\|_{(W^{1,p}(\Omega;\R^3))'}$.
In particular, by the definition of $\mathcal{Y}$ we infer that
$\{y_k\}$ is a bounded sequence in $W^{1,p}(\Omega; \R^3)$, so that
there exists $y \in \mathcal{Y}$ such that
 \begin{equation}y_k\wk y \quad 
 \text{ in } W^{1,p}(\Omega; \R^3)\label{eq:prima}
\end{equation}
 for some not relabeled subsequence. Since $G^{-1}\in
 L^\infty(\Omega;\R^{3\times 3})$, we deduce that
 \begin{equation}
 \label{eq:conv-grady-G0}
 \nabla y_k G ^{-1}\wk \nabla y \,G ^{-1}\quad 
\text{ in }  L^p(\Omega; \R^{3\times 3}).
 \end{equation}
 The weak continuity of the minors of $\nabla y_k$ \cite{ball} for $p>3$ and \eqref{eq:conv-grady-G0} 
yield
\begin{align*}
&{\rm{cof }} ( \nabla y_k \, G ^{-1}) \wk  {\rm{cof }} ( \nabla y 
                 \,G ^{-1}) \quad 
\text{ in } L^{p/2}(\Omega; \R^{3\times 3}), \\
&\det ( \nabla y_k \,G ^{-1}) \wk  \det( \nabla y  \,G ^{-1}) \quad 
\text{ in } L^{p/3}(\Omega),
 \end{align*} 
which, combined with the polyconvexity from {\bf (H1)}, imply 
 \begin{equation}
 \label{eq:W}
\int_{\Omega} W(\nabla y  \, G^{-1})\,  \det G   \,\di x \leq
\liminf_{k\to \infty} \int_{\Omega} W(\nabla y_k \, G ^{-1})\, 
\det G  \,\di x.
 \end{equation}
By  \eqref{eq:prima} we have $\langle\ell ,y_k\rangle \to \langle\ell
,y \rangle$, 
which, together with \eqref{eq:W}, leads to the minimality \eqref{eq:elastic}.
\end{proof}\EEE

\RRR Owing to Lemma \ref{lemma:elastic}, for every $G^0\in
\mathcal{G}_\infty$  with $\det G^0 \geq \delta>0$ a.e.  there exists
$y^0\in \mathcal{Y}$ solving \eqref{eq:discrete-growth1} for $i=0$. \EEE
The existence of discrete solutions is \RRR then \EEE guaranteed by the following lemma.

\begin{lemma}[Discrete existence]
  \label{lemma:exist-y}
  \RRR Under assumptions {\bf(H1)}--{\bf(H3)}, {\bf(H5)}--{\bf(H6)},
  let $(y_0,G_0) = (y^0,G^0)$, where $y^0$ solves \eqref{eq:elastic}
  for $G=G^0$ and $\ell=\ell_0$. \EEE For 
   every $i=1,\dots, N$, there exists
 $(y_i,G_i)\in \mathcal{Y}\times \mathcal{G}_\infty$ solving
 \eqref{eq:discrete-growth1}--\eqref{eq:discrete-growth2} \UUU with
 $\det G_i \geq \exp(-3\tau i \| M \|_{L^\infty}) \delta$ a.e. \EEE
 \end{lemma}

 \begin{proof}  \RRR
   We proceed by induction on $i$. Assume that there exist $G_{i-1}\in \mathcal{G}_\infty$ with
   $\det G_{i-1} \geq \exp(-3\tau(i-1) \| M \|_{L^\infty}) \delta$
   a.e. and define $G_i$ via \EEE position  \eqref{eq:discrete-growth2}.   Let us \EEE check that $G_i\in W^{1,\infty}(\Omega; \R^{3\times 3})$, and that
 $\det G_i \UUU \geq  
\exp(-3\tau i \| M \|_{L^\infty}) \delta $ \EEE  a.e. in $\Omega$, for $i=1,\dotsc,N$. 
 By differentiating \eqref{eq:discrete-growth2}, we find 
 \begin{equation}
 \label{eq:nabla-G}
 \nabla G_i=(G_{i-1}^\top \nabla  \exp \big( \tau \,  M(G_{i-1}, (K_\tau {\nabla y})_{i-1}) \big)^\top)^t
 +\exp \big(\tau \,  M(G_{i-1}, (K_\tau {\nabla y})_{i-1}) \big)\nabla G_{i-1}
\end{equation}
a.e. in $\Omega$, for $i=1,\dotsc,N$, where the above equality holds in the sense of distributions.
 In view of {\bf (H6)} and by the fact that $G_{i-1}\in \mathcal{G}_\infty$, we infer the estimate
 
 \begin{align*}
\|\nabla G_i \|_{L^\infty}
 & \leq  \| G_{i-1}^\top \nabla \exp \big( \tau \,  M(G_{i-1}, (K_\tau {\nabla y})_{i-1}) \big)^\top \|_{L^\infty}  \\ 
 & \quad +  \| \exp\big(\tau \,  M(G_{i-1}, (K_\tau {\nabla y})_{i-1}) \big)\nabla G_{i-1} \|_{L^\infty}  \\[0.2cm]
 & \leq  \| G_{i-1}\|_{L^\infty} \| \nabla \exp \big( \tau \,  M(G_{i-1}, (K_\tau {\nabla y})_{i-1}) \big)\|_{L^\infty}   \\ 
 & \quad +  \| \exp\big(\tau \,  M(G_{i-1}, (K_\tau {\nabla y})_{i-1}) \big) \|_{L^\infty}
  \| \nabla G_{i-1} \|_{L^\infty}  \\[0.2cm]
 &\leq C\left(1+ \|\nabla  \exp \big(\tau \,  M(G_{i-1}, (K_\tau {\nabla y})_{i-1})\|_{L^\infty}\right)
 \end{align*} 
with $C$ depending on $\| M \|_{L^\infty}$ and $\| G_{i-1} \|_{W^{1,\infty}}$.
Denoting by $D_1 M$ and $D_2 M$ the differentials of $M$ with respect to its first and second matrix-valued variables, respectively, we obtain from 
the properties of the matrix exponential that 
\begin{align}
&\nonumber|\nabla  \exp \big(\tau \,  M(G_{i-1}, (K_\tau {\nabla y})_{i-1})\big) | 
\\[0.2cm]
&\nonumber \quad = 
\tau\bigg| \int_0^1 \exp \Big( (1-s) \tau \, M \big(G_{i-1}, (K_\tau {\nabla y})_{i-1}\big)\Big) \, \nabla M (G_{i-1}, (K_\tau {\nabla y})_{i-1}) \, 
\\
&\nonumber \hspace{1.5cm}  \times \exp \Big( s \tau \, M \big(G_{i-1}, (K_\tau {\nabla y})_{i-1}\big)\Big) \, \di s\bigg| 
\\[0.2cm]
&\nonumber\quad \leq \tau\int_0^1 \exp \Big( (1-s) \tau \big|M \big (G_{i-1}, (K_\tau {\nabla y})_{i-1}\big) \big| \Big) \, 
\exp \Big( s \tau \big|M \big (G_{i-1}, (K_\tau {\nabla y})_{i-1}\big) \big| \Big) \, \di s \, 
\\
&\nonumber  \hspace{1.5cm} \times \big|\nabla M \big (G_{i-1}, (K_\tau {\nabla y})_{i-1}\big) \big| 
\\[0.2cm]
&\nonumber\quad = \tau\exp \Big(  \tau \big|M \big (G_{i-1}, (K_\tau {\nabla y})_{i-1}\big) \big| \Big) \, 
\Big| D_1 M \big(G_{i-1}, (K_\tau {\nabla y})_{i-1}\big) \, \nabla G_{i-1} 
\\
&\nonumber \hspace{6cm} +
D_2 M \big (G_{i-1}, (K_\tau {\nabla y})_{i-1}\big) \, \nabla (K_\tau {\nabla y})_{i-1} \Big| 
\\
&\label{eq:est-grad-exp}\quad \leq 2\tau \exp \big(  \tau\Vert M\Vert_{L^{\infty}} \big) 
\Vert M\Vert_{W^{1,\infty}} 
\Big( |\nabla G_{i-1} | + |\nabla (K_\tau {\nabla y})_{i-1}| \Big),
\end{align} 
so that, using once again {\bf (H6)} and the fact that
$G_{i-1}\in \mathcal{G}_\infty$  we have that 
$$
\|\nabla G_i \|_{L^\infty} 
\leq C \tau (1+  \|\nabla (K _\tau {\nabla y})_{i-1}\|_{L^\infty}) +C, 
$$
where $C$ depends on $\Vert M\Vert_{W^{1,\infty}}$, and $\| G_{i-1} \|_{W^{1,\infty}}$. 
Now, by \eqref{eq:discrete-conv} and {\bf (H4)},
\begin{equation}
\label{eq:est-grad-conv}
\|\nabla (K_\tau {\nabla y})_{i-1}\|_{L^\infty}\leq \sum_{j=0}^{i-1}\tau  | \kappa_j | 
 \|{\nabla y}_{i-1-j}\|_{L^p(\Omega; \R^{3\times 3})}\|\nabla \phi\|_{L^q(\R^3; \R^{3})}\leq C
\end{equation}
owing to the fact that $y_{j}\in \mathcal{Y}$ for $j=0,\dots,i-1$. This yields that $G_i\in W^{1,\infty}(\Omega; \R^{3\times 3})$.

\UUU The lower bound on the \EEE determinant of $G_i$ follows by induction, namely, 
\begin{align}
	\det G_i &= \exp \big( \tau \, {\rm tr} \, M(G_{i-1}, (K_\tau {\nabla y})_{i-1})\big) \, 
	\det G_{i-1} \label{eq:bd-det-below} \geq \exp \big( -3\tau \,
                            \Vert M\Vert_{L^{\infty}} \big) \, \det
                   \, G_{i-1} \\
  &\UUU \geq \exp \big( -3\tau \,
                            \Vert M\Vert_{L^{\infty}} \big)
                  \exp \big( -3\tau (i-1) \,
                            \Vert M\Vert_{L^{\infty}} \big) \delta  =\exp \big( -3\tau i \,
                            \Vert M\Vert_{L^{\infty}} \big) \delta
\end{align}
a.e. in $\Omega$, for $i=1,\dotsc,N$.

\RRR
We can hence conclude the proof by applying Lemma \ref{lemma:elastic} for $G=G_{i}$ and
$\ell=\ell_{i}$ in order to find a deformation $y_i\in \mathcal{Y}$
\EEE solving~\eqref{eq:discrete-growth1}. 
\end{proof}

\subsection{A-priori estimates}
\label{subs:estimates}
Denoting by $\overline{(K_\tau {\nabla y})}_\tau$ and
$\widehat{(K_\tau {\nabla y})}_\tau$ the backward piecewise
constant and piecewise affine interpolants associated to $K_\tau
{\nabla y}$, cf. \eqref{eq:pw-aff}, \eqref{eq:pw-con}, and \eqref{eq:discrete-conv}, the main result of this subsection is the following.
\begin{proposition}[A-priori estimates]
	\label{prop:apriori}
There exist $\tau^*\in (0,1)$ depending on $\Vert M\Vert_{L^{\infty}}$  such that, for every
$\tau\in (0,\tau^*)$ we have 
\begin{align}
&\label{eq:hatG-winfty}
	\|\widehat{G}_{\tau}\|_{W^{1,\infty} (0,T; L^{\infty} (\Omega; \R^{3\times 3}))}\leq C,\\ 
&\label{eq:unif-bd-det}
{\rm det}\,\overline{G}_\tau(t)\geq C{\rm det}\,G_0 \UUU \geq C\delta
                                                                                                     \quad
                                                                                                     \text{a.e. in}
                                                                                                     \
                                                                                                     \Omega,
                                                                                                     \
                                                                                                     \EEE
                                                                                                     \forall
                                                                                                     \, t\in [0,T],\\
&\label{eq:unif-y} \|\overline{ y}_{\tau}\|_{L^\infty(0,T; W^{1,p}(\Omega; \R^{3}))}+\|\widehat{y}_\tau\|_{L^\infty(0,T;W^{1,p}(\Omega; \R^{3}))}\leq C,\\
&\label{eq:est-nabla-g-linfty}
\|\overline{ G}_{\tau}\|_{L^\infty(0,T;W^{1,\infty}(\Omega; \R^{3\times 3}))}+\|\widehat{ G}_\tau\|_{L^\infty(0,T; W^{1,\infty}(\Omega; \R^{3\times 3}))}\leq C,\\
&\label{eq:bd-disc-conv1} 
\|\overline{(K_\tau {\nabla y})}_\tau\|_{L^\infty(0,T;W^{1,\infty}(\Omega; \R^{3\times 3}))}\leq C,  \\
&\label{eq:bd-disc-conv2} 
\|\widehat{(K_\tau {\nabla y})}_\tau\|_{L^\infty(0,T;W^{1,\infty}(\Omega; \R^{3\times 3}))}+\|\widehat{(K_\tau {\nabla y})}_\tau\|_{W^{1,\infty}(0,T;L^{\infty}(\Omega;\ \R^{3\times 3}))}\leq C, 
\end{align} 
where the positive constant $C$ depends on  $c_1$, $\Omega$,  $\Vert G^0
\Vert_{W^{1,\infty}}$, $\Vert M\Vert_{W^{1,\infty}}$,  $\|
\kappa\|_{W^{1,1}(0,T)}$, $\| \phi\|_{W^{1,q}(\R^3)}$, $\|
\ell\|_{W^{1,1}(0,T;(W^{1,p}(\Omega;\R^3))')}$, and~$T$. 
\end{proposition}
\begin{proof}
Within this proof, the symbol $C$ stands for a positive
constant, possibly depending on  $c_1$, $\Omega$, $\Vert G^0
\Vert_{W^{1,\infty}}$, $\Vert M\Vert_{W^{1,\infty}}$,  $\|
\kappa\|_{W^{1,1}(0,T)}$, $\| \phi\|_{W^{1,q}(\R^3)}$, $\|
\ell\|_{W^{1,1}(0,T;(W^{1,p}(\Omega;\R^3))')}$, and  $T$ but
independent of $\tau$. The actual value of $C$ can change from line to
line.

We first show estimate \eqref{eq:hatG-winfty}. We subtract $G_{i-1}$ from both sides of \eqref{eq:discrete-growth2}, divide by $\tau$, and contract with $G_i$. This yields
\begin{equation}
\label{eq:diff-quot}
\frac{G_i - G_{i-1}}{\tau} : G_i = \frac{ \Big( \exp \big( \tau \,  M(G_{i-1}, (K_\tau {\nabla y})_{i-1})\big) - {\rm Id} \Big)  G_{i-1}}{\tau} : G_i 
\end{equation}
a.e. in $\Omega$, for $i=1,\dotsc,N$. 
The left-hand side can be treated as 
\begin{equation*}
\frac{G_i - G_{i-1}}{\tau} : G_i = \frac{ |G_i|^2}{2\tau} + \frac{ |G_i- G_{i-1}|^2}{2\tau} 
- \frac{ |G_{i-1}|^2}{2\tau}, 
\end{equation*}
whereas for the right-hand side of \eqref{eq:diff-quot} we use
the Cauchy-Schwarz and the Young inequalities, as well as the 
properties of the matrix exponential in order to get that 
\begin{align*}
&	\frac{ \Big( \exp \big( \tau \, M(G_{i-1}, (K_\tau {\nabla y})_{i-1})\big) - {\rm Id} \Big)  G_{i-1}}{\tau} : G_i  
\\
& \hspace{1.5cm} \leq 
    \left|  \frac{ \exp \big( \tau \,  M(G_{i-1}, (K_\tau {\nabla y})_{i-1})\big) - {\rm Id} }{\tau}\right| \, | G_{i-1}| \, |G_i|
 \\
&\hspace{1.5cm} \leq \frac{ \exp \big( \tau \, \left| M(G_{i-1}, (K_\tau {\nabla y})_{i-1}) \right| \big) -  1}{\tau} \, | G_{i-1}| \, |G_i| 
 \\
&\hspace{1.5cm}\leq  \frac{ \exp \big( \tau \,  \Vert M\Vert_{L^{\infty} } \big) -  1}{\tau} \left( \frac{| G_{i-1}|^2 }{2} + \frac{| G_{i}|^2 }{2} \right)
\end{align*}
 where we have also used the fact that $|\exp A -  {\rm Id}| \leq
\exp(|A|)-1$. 
Therefore, we have obtained that
\begin{equation*}
\frac{ |G_i|^2}{2\tau} + \frac{ |G_i- G_{i-1}|^2}{2\tau} 
- \frac{ |G_{i-1}|^2}{2\tau} \leq \frac{ \exp \big( \tau \,  \Vert M\Vert_{L^{\infty}} \big) -  1}{\tau} \left( \frac{| G_{i-1}|^2 }{2} + \frac{| G_{i}|^2 }{2} \right).
\end{equation*}
Fix an integer $m \leq N$.  By multiplying by $\tau$ and summing up for $i=1,\dotsc,m$, we deduce that 
\begin{align*}
\frac{ |G_m|^2}{2}- \frac{ |G_0|^2}{2} 
&\leq \frac{1}{2} \left( \exp \big( \tau \,  \Vert M\Vert_{L^{\infty}
  } \big) -1 \right) \left(  |G_m|^2 + |G_0|^2 \right) 
 + \sum_{i=1}^{m-1} \left( \exp \big( \tau \,  \Vert M\Vert_{L^{\infty}} \big) -1 \right) |G_i|^2 . 
\end{align*}
Taking $\tau<(\log 2)/\Vert M\Vert_{L^{\infty}  } =:\tau^*$ and applying the Discrete Gronwall Lemma, cf. \cite[Proposition 2.2.1]{jerome}, we conclude that 
\begin{equation}
	\label{eq:Gbound}
|G_m| \leq C, \quad \forall \, m=1,\dotsc,N, 
\end{equation}
and hence 
\begin{equation}
\label{eq:l-infty-g}
	\widehat{G}_{\tau} \text{ and } \overline{G}_{\tau} \text{ are
          bounded in }L^{\infty} (0,T; L^{\infty} (\Omega; \R^{3\times
          3})) \text{ independently of } \tau  \in (0, \tau^*).
\end{equation}

By subtracting $G_{i-1}$ from both sides of \eqref{eq:discrete-growth2}, dividing by $\tau$, and taking the norm, we get 
\begin{align}
\nonumber\left| \frac{G_i - G_{i-1}}{\tau} \right| 
&\leq \left| \frac{ \exp \big( \tau \, M(G_{i-1}, (K_\tau {\nabla y})_{i-1})\big) - {\rm Id}}{\tau} \right| \, |G_{i-1}|
\\
&\label{eq:ing0} \leq \frac{1}{\tau}  \left( \exp \big( \tau \,  \Vert M\Vert_{L^{\infty}} \big) -1 \right) |G_{i-1}| 
\quad \forall \, i=1,\dotsc, N. 
\end{align}
By \eqref{eq:Gbound}, for $\tau<  \tau^\ast $ we infer \eqref{eq:hatG-winfty}.

 Next, we prove \eqref{eq:unif-bd-det} and \eqref{eq:unif-y}. 
Recalling the definition of $\mathcal{G}_\infty$, by subtracting
$G_{i-1}$  from both sides of \eqref{eq:discrete-growth2}, dividing by $\tau$, taking the gradient, and contracting with $\nabla G_i$, we obtain 
\begin{equation}
\label{eq:ing1}
\nabla \left( \frac{G_i - G_{i-1}}{\tau} \right) : \nabla G_i = \nabla
( E_{i-1}  G_{i-1}) : \nabla G_i 
\end{equation} 
a.e. in $\Omega$, for $i=1,\dotsc,N$, 
where we have set 
\begin{equation}
\label{eq:ing2}
	E_{i-1} := \frac{ \exp \big( \tau \, M(G_{i-1}, (K_\tau {\nabla y})_{i-1})\big) - {\rm Id}}{\tau}. 
\end{equation}
 The left-hand side of \eqref{eq:ing1} reads equivalently 
\begin{equation}
\label{eq:ing3}
\nabla \left( \frac{G_i - G_{i-1}}{\tau} \right) : \nabla G_i = \frac{ |\nabla G_i|^2}{2\tau} + \frac{ |\nabla G_i- \nabla G_{i-1}|^2}{2\tau} - \frac{ |\nabla G_{i-1}|^2}{2\tau}, 
\end{equation}
while for the right-hand side, we compute 
\begin{align}
\nonumber	\nabla (E_{i-1} G_{i-1}) : \nabla G_i 
&\leq | \nabla( E_{i-1}G_{i-1} )|  |\nabla G_i| 
=  | (G_{i-1}^\top \nabla E_{i-1}^\top)^t  + E_{i-1}\nabla G_{i-1} |  |\nabla G_i| 
\\[0.1cm]
&\nonumber \leq \big( |G_{i-1}^\top \nabla E_{i-1}^\top| + |E_{i-1}\nabla G_{i-1} | \big)  |\nabla G_i| 
\\[0.1cm]
&\label{eq:ing4} \leq \big( |G_{i-1}| | \nabla E_{i-1}| + |E_{i-1}| | \nabla G_{i-1} | \big)  |\nabla G_i| .
\end{align}
From the properties of the matrix exponential, arguing as in \eqref{eq:est-grad-exp}, 
 and from the estimate on the convolution \eqref{eq:est-grad-conv} we deduce that 
\begin{align}
\nonumber&|\nabla E_{i-1}| 
\leq 2 \exp \big(  \tau \|M\|_{L^{\infty}} \big) 
\Vert M\Vert_{W^{1,\infty}} 
\Big( |\nabla G_{i-1} | + |\nabla  (K_\tau {\nabla y})_{i-1}| \Big)\\
& \leq 2 \exp\big(\tau \|M\|_{L^{\infty}})\Vert
  M\Vert_{W^{1,\infty}}
\Big( |\nabla G_{i-1} | +
  \sum_{j=0}^{i-1}\tau |\kappa_j |  \|{\nabla y}_{i-1-j}\|_{L^p(\Omega;\R^{3\times 3})}\|\nabla \phi\|_{L^q(\R^3;\R^{3})}\Big). \label{eq:nabla-ni}
\end{align}
On the other hand, by iterating \eqref{eq:bd-det-below}, we find
$$
{\rm det}\,G_{ m}\geq  {\rm exp} (-3\tau m \| M
\|_{L^\infty}) \, {\rm det}\,G_0 \geq  {\rm exp} (-3 T \| M
\|_{L^\infty}) \, {\rm det}\,G_0\quad \forall \, m =1,\dots,N, 
$$
which gives \eqref{eq:unif-bd-det}. 
Thus,  {\bf (H2)},  {\bf (H6)},  \eqref{eq:discrete-growth1}, and the Poincar\'e inequality imply
\begin{align*}
\|\nabla y_{i}\, G_i^{-1}\|^p_{L^p(\Omega;\R^{3\times 3})}&\leq C\left( 1 + \langle \ell_i, y_i \rangle \right) \leq C\left( 1 +  \| \ell_i \|_{(W^{1,p} (\Omega; \R^3))'}  
\| y_i \|_{W^{1,p} (\Omega; \R^3)}  \right) 
\\
& \leq C \left( 1+ \|\nabla y_{i}\|_{L^p(\Omega; \R^{3\times 3})} \right) 
\\
& \leq  C \left( 1+ \|\nabla y_{i}\,G_i^{-1}\|_{L^p(\Omega;\R^{3\times 3})} 
\| G_i \|_{L^{\infty}} \right) . 
\end{align*}
By using the Young Inequality and \eqref{eq:l-infty-g}, we arrive at 
\begin{equation}
\label{eq:new-est-strain}
\|\nabla y_{i}G_i^{-1}\|^p_{L^p(\Omega;\R^{3\times 3})} \leq C \left( 1+ \| G_i \|^q_{L^{\infty}}  \right) \leq C,
\end{equation}
which in turn yields
\begin{equation}
	\label{eq:est-y-tau}
\|\nabla y_{i}\|^p_{L^p(\Omega;\R^{3\times 3})}\leq  \|\nabla y_{i}\,G_i^{-1}\|^p_{L^p(\Omega;\R^{3\times 3})} \|G_i\|^p_{L^{\infty}} \leq C \quad \forall \, i=1,\dotsc, N .
\end{equation}
 In particular, we obtain the bound \eqref{eq:unif-y}.
 
 We now prove the bound \eqref{eq:est-nabla-g-linfty}. 
Going back to \eqref{eq:nabla-ni}, from \eqref{eq:unif-y} and  {\bf (H4)} we infer the estimate
$$
|\nabla E_{i-1}|\leq C(1+|\nabla G_{i-1}|).
$$
 On the other hand, the same computations as in \eqref{eq:ing0} yield
$$
|E_{i-1}|\leq \frac{1}{\tau}\left(\exp(\tau \|M\|_{L^\infty})-1\right)
\leq C. 
$$
Therefore, by combining \eqref{eq:ing1}--\eqref{eq:ing4} and
multiplying by $\tau$,  from  \eqref{eq:l-infty-g} we conclude that
$$
\frac{ |\nabla G_i|^2}{2} + \frac{ |\nabla G_i- \nabla G_{i-1}|^2}{2}
- \frac{ |\nabla G_{i-1}|^2}{2}\leq C\tau \left( 1+ |\nabla
  G_{i-1}| \right) |\nabla G_i|  
$$
for $i=1,\dotsc,N$. Summing up for $i=1,\dotsc,m$ for some
$m=1,\dotsc,N$ and using the Young Inequality we get 
$$
\frac{ |\nabla G_m|^2}{2} - \frac{ |\nabla G_0|^2}{2} 
\leq  \sum_{i=1}^m C\tau (1+ |\nabla G_{i-1}|)|\nabla G_i| \leq
\sum_{i=1}^m C\tau (1+ |\nabla G_{i-1}|^2) + \frac14 \sum_{i=1}^m\tau
|\nabla G_i|^2. 
$$ 
We can hence apply the  Discrete Gronwall Lemma and obtain $|\nabla G_m| 
\leq C$ for all $m=1,\dotsc,N$, which, together with \eqref{eq:l-infty-g}, implies  \eqref{eq:est-nabla-g-linfty}.

It remains to prove the a-priori bounds involving the discrete convolution term. From \eqref{eq:est-grad-conv}, \eqref{eq:unif-y}, 
and {\bf (H4)}, for every $i=1,\dots,N$ there holds
\begin{align}
&\nonumber \|\nabla (K_\tau {\nabla y})_i\|_{L^\infty}
\leq \sum_{j=0}^N \tau |\kappa_j| \|\nabla y_{i-j}\|_{L^p(\Omega; \R^{3\times 3})}\|\nabla \phi\|_{L^q(\R^3; \R^{3})}\\
&\quad \label{eq:est-grad-disc-conv1}
\leq \tau (N+1) \|\kappa\|_{L^\infty(0,T)}\|\nabla \overline{y}_\tau\|_{L^{\infty}(0,T;L^p(\Omega;\R^{3\times 3}))}\|\nabla \phi\|_{L^q(\R^3; \R^{3})}\leq C.
\end{align} 
Analogously,
\begin{align*}
\frac{(K_\tau {\nabla y})_i(x)-(K_\tau {\nabla y})_{i-1}(x)}{\tau}=\int_{\R^3}\phi(x-z)\left( \kappa_0{\nabla y}_i(z)+
\sum_{j=1}^i (\kappa_j-\kappa_{j-1}){\nabla y}_{i-j}(z)\right) \di z
\end{align*}
for a.e. $x\in \Omega$, so that
\begin{align*}
\left|\frac{(K_\tau {\nabla y})_i-(K_\tau {\nabla y})_{i-1}}{\tau}\right|\leq (\|\kappa\|_{L^\infty(0,T)}+\|\kappa'\|_{L^1(0,T)})\|\nabla \overline{y}_\tau\|_{L^\infty(0,T;L^p(\Omega;\R^{3\times 3}))}\|\phi\|_{L^q(\R^3)}
\end{align*}
a.e. in $\Omega$. This, combined with \eqref{eq:unif-y} and  {\bf (H4)}, yields \eqref{eq:bd-disc-conv1} and \eqref{eq:bd-disc-conv2} and concludes the proof of the proposition.
\end{proof}

\subsection{Passage to the limit}
\label{subs:limit}
We are now in a position to  pass to the limit as $\tau \to
0$. 
In view of \eqref{eq:hatG-winfty}, \eqref{eq:unif-y}, and \eqref{eq:est-nabla-g-linfty},
we  find  $\allowbreak G\in W^{1,\infty}(0,T;L^\infty(\Omega;\R^{3\times 3}))\cap L^\infty(0,T;W^{1,\infty}(\Omega;\R^{3\times 3}))$ and a limiting deformation $y\in L^{\infty}(0,T;W^{1,p}(\Omega;\R^3))$ such that, up to the extraction of not relabeled subsequences, there holds
\begin{align}
\label{eq:conv1}
\widehat{G}_\tau\wstarto G\quad
&\text{ in } W^{1,\infty}(0,T;L^\infty(\Omega;\R^{3\times 3}))\cap L^\infty(0,T;W^{1,\infty}(\Omega;\R^{3\times 3})),\\
\label{eq:conv1bis}\overline{y}_\tau\wstarto y\quad 
&\text{ in } L^{\infty}(0,T;W^{1,p}(\Omega;\R^3)).
\end{align}
 The Aubin-Lions Lemma \cite{Simon87} yields
\begin{equation}
\label{eq:conv2}
\widehat{G}_\tau\to G\quad
\text{ in } C([0,T];C^\alpha(\overline{\Omega};\R^{3\times 3}))\quad \forall \, \alpha\in (0,1),
\end{equation}
and, in particular,
\begin{equation}
\label{eq:in-cond-G}
G(0)=G_0.
\end{equation}
From \eqref{eq:pw-aff} and \eqref{eq:pw-con}, we deduce the identity
\begin{align*}
|\overline{G}_\tau(t)-\widehat{G}_\tau(t)|=\tau(1-\alpha_i(t))|\widehat{G}_\tau'(t)|\quad\text{ a.e. in } \Omega,
\end{align*}
for all $t \in ((i-1)\tau,i\tau]$, $i=1,\dots,N$, so that
\begin{equation}
\label{eq:Gtaubound}
|\overline{G}_\tau(t)-\widehat{G}_\tau(t)|\leq \tau |\widehat{G}_\tau'(t)|    \quad\text{ a.e. in } \Omega, \,
\forall \, t\in [0,T].
\end{equation}
Hence, \eqref{eq:hatG-winfty}, \eqref{eq:est-nabla-g-linfty}, and \eqref{eq:conv2} yield
\begin{align}
&\label{eq:conv3} \overline{G}_\tau\to G\quad
                \text{ in }  L^\infty ((0,T) \times \Omega ;\R^{3\times 3}),\\
  &\label{eq:conv3b}  \overline{G}_\tau(t)\to G(t)\quad
 \text{ in }  L^\infty (\Omega ;\R^{3\times 3}), \ \forall t \in[0,T],\\
&\label{eq:conv}\overline{G}_\tau\wstarto G\quad
\text{ in } L^\infty(0,T;W^{1,\infty}(\Omega;\R^{3\times 3})).
\end{align}
Moreover,  by \eqref{eq:unif-y}, we find $\xi \in L^\infty(0,T;L^p(\R^3 ;\R^{3\times 3}))$ such that, up to subsequences,
$$
{\nabla \overline{y}_\tau}\wstarto \xi \quad 
\text{ in } L^{\infty}(0,T;L^{p}(\R^3;\R^{3\times 3}))
$$
where ${\nabla \overline{y}_\tau}$ is the trivial extension to
the whole $\R^3$. 

We proceed by identifying $\xi$. Let $x_0\in \Omega$ and $r>0$ be such
that $B_r(x_0) := \{x \in \R^3 \ : \ |x-x_0|<r\}  \subset\Omega$. Then, for every $\psi\in L^1((0,T);\R^{3\times 3})$ we have
$$
\int_0^T\int_{B_r(x_0)} \xi(t,x)\cdot \psi(t)\,\di{x}\,\di{t}=\lim_{\tau\to 0}\int_0^T\int_{B_r(x_0)} \nabla \overline{y}_\tau(t,x)\cdot \psi(t)\,\di{x}\,\di{t}=\int_0^T\int_{B_r(x_0)} \nabla y(t,x)\cdot \psi(t)\,\di{x}\,\di{t},
$$
where the last equality is due to \eqref{eq:conv1bis}. 
This yields $\int_{B_r(x_0)}\xi(t,x)\,\di{x}=\int_{B_r(x_0)}\nabla y(t,x)\,\di{x}$ for a.e. $t\in (0,T)$ and for every $B_r(x_0)\subset \Omega$. Namely, $\xi=\nabla y$ a.e. in $(0,T)\times \Omega$. Analogously, let $y_0\in \Omega^c$ and let $s>0$ be such that $B_s(y_0)\subset \Omega^c$. Then, for every $\psi\in L^1((0,T);\R^{3\times 3})$ there holds
$$
\int_0^T\int_{B_s(y_0)} \xi(t,x)\cdot \psi(t)\,\di{x}\,\di{t}=\lim_{\tau\to 0}\int_0^T\int_{B_s(y_0)} \nabla \overline{y}_\tau(t,x)\cdot \psi(t)\,\di{x}\,\di{t}=0,
$$
which implies $\xi=0$ a.e. in $(0,T)\times \Omega^c$. Therefore, $\xi={\nabla y}$.

Arguing as in \eqref{eq:diff-quot}, we deduce
\begin{equation}
\label{eq:diff-quot2}
\widehat{G}_\tau'(t) = \frac{ \Big( \exp \big( \tau \, M(\overline{G}_\tau(t-\tau), \overline{(K_\tau {\nabla y})}_{\tau}(t-\tau))\big) - {\rm Id} \Big) \, }{\tau} \, \overline{G}_{\tau}(t-\tau).
\end{equation}
Owing to \eqref{eq:bd-disc-conv1} and \eqref{eq:bd-disc-conv2}, the same argument as in the proof of \eqref{eq:conv3} yields
\begin{equation}
\label{eq:claim-disc-conv}
\overline{(K_\tau {\nabla y})}_{\tau}\to {K {\nabla y}}\quad
\text{ in }   L^\infty((0,T)\times  \Omega;\R^{3\times 3}).
\end{equation}
Therefore, by \eqref{eq:conv1}, \eqref{eq:conv3}, the regularity of
$M$, and the fact that the exponential map is locally
Lipschitz, we find that $(y,G)$ solves \eqref{eq:problemG}. 

 We are hence left with proving that minimality
\eqref{eq:problemP} holds \RRR almost everywhere in time. \EEE To this
aim, assume to be given $\hat y\in
\mathcal{Y}$ and \RRR recall that $ W(\nabla \hat y \, G^{-1}(t)) \
\det G(t) \in L^\infty(\Omega)  $ since $G(t) \in
L^\infty(\Omega ; \R^{3\times 3})$ and $\det G(t) \geq C\delta$ a.e. In particular, we have that
\begin{equation}
  \label{eq:controllo}
  W(\nabla \hat y \, G^{-1}(t)) \
\leq \MMM \lambda \RRR(\| G(t) \|_{L^{\infty}}, C\delta)\quad
\text{a.e. in} \ \Omega \times
(0,T)
\end{equation}
where \MMM $\lambda$ \RRR is defined in \eqref{k}. \EEE We will make use
of the following 
simplified version of \cite[Lemma 4.1]{ms5}.
\begin{lemma}\label{lemma:ball}
  Under assumption {\bf(H2)} there exist $c_3, \, \epsilon>0$ such that
  \begin{equation}
    |W(FH)- W(F)| \leq c_3 (W(F)+1) | H- {\rm Id}| \quad \forall F,
    \, H \in GL_+(3), \ |H-{\rm Id}|\leq \epsilon.\label{eq:ball}
  \end{equation}
\end{lemma}
\RRR Let $s\in (0,T)$. \EEE By choosing $F=\nabla \hat y \, G^{-1}(s)$ and $H = G(s)
\overline{G}_\tau^{-1}(s)$ in \eqref{eq:ball} we find that
\begin{align}
& \int_\Omega W(\nabla \hat y\,\overline{G}_\tau^{-1}(s))\,\UUU \det \overline{G}_\tau(s)\EEE\, \di{x}  -
  \int_\Omega W(\nabla \hat y\,{G}^{-1}(s))\, \UUU \det G(s)\EEE\,\di{x} \nonumber \\
  &\quad \leq
  \left(c_3\int_\Omega(W(\nabla \hat y\,{G}^{-1}(s))+1) \, \UUU \det \overline{G}_\tau(s)\EEE\,\di{x}\right) \|
    G(s) \overline{G}_\tau^{-1}(s) - {\rm Id} \|_{L^\infty} \nonumber\\
  &\UUU \qquad + \int_\Omega W(\nabla \hat y\,{G}^{-1}(s))\, \UUU \big(\det \overline{G}_\tau(s)-\det G(s)\big)\EEE\,\di{x} . \nonumber
\end{align}
\RRR Owing to \eqref{eq:controllo}, \EEE this entails that $W(\nabla \hat y \overline{G}_\tau^{-1}(t)) \, \RRR
\det  \overline{G}_\tau(t) \EEE \in
L^1(\Omega)$ for all $\tau $ as well. Moreover, taking into account convergence
\eqref{eq:conv3b}, we have that 
\begin{align}
  \label{eq:FH1}
&\limsup_{\tau \to 0}\left(\int_\Omega W(\nabla \hat y\,\overline{G}_\tau^{-1}(s))\,\UUU \det \overline{G}_\tau(s)\EEE\, \di{x}  -
  \int_\Omega W(\nabla \hat y\,{G}^{-1}(s))\, \UUU \det G(s)\EEE\,\di{x} \right) =0.
\end{align}

\RRR Fix now $t\in (0,T)$ and $\delta>0$ small so that $t+\delta \in
(0,T)$. 
\EEE  From the discrete minimality  \eqref{eq:discrete-growth1}, we
  deduce that 
\begin{equation}
\RRR \int_t^{t+\delta } \EEE \!\! \int_{\Omega}W(\nabla
\overline{y}_\tau(s)\,\overline{G}_\tau^{-1}(s))\, \UUU \det
\overline{G}_\tau(s)\EEE\,\di{x}\, \RRR\di{s}\EEE-\RRR \int_t^{t+\delta } \EEE \!\! \langle \overline{\ell}_\tau(s),\overline{y}_\tau(s)\rangle \, \RRR\di{s}\EEE\leq \RRR \int_t^{t+\delta } \EEE \!\! \int_{\Omega}W(\nabla \hat y\,\overline{G}_\tau^{-1}(s))\, \UUU \det \overline{G}_\tau(s)\EEE\,\di{x}\, \RRR\di{s}\EEE-\RRR \int_t^{t+\delta } \EEE \!\! \langle \overline{\ell}_\tau(s),\hat y\rangle \, \RRR\di{s}\EEE.\label{eq:dispo}
\end{equation}
\RRR We now aim at passing to the $\liminf$ as $\tau \to 0$ first. \EEE
Convergence \eqref{eq:conv1bis} and the regularity of the
applied loads implies that
$$\RRR\int_t^{t+\delta } \EEE \!\! \langle \overline{\ell}_\tau(s),\overline{y}_\tau(s)\rangle \, \RRR\di{s}\EEE\to
\RRR\int_t^{t+\delta } \EEE \!\! \langle {\ell}(s),{y}(s)\rangle \,
\RRR\di{s}\EEE\quad \text{and}\quad \RRR\int_t^{t+\delta } \EEE \!\!
\langle \overline{\ell}_\tau(s),\hat y\rangle \, \RRR\di{s}\EEE\to
\RRR\int_t^{t+\delta } \EEE \!\! \langle {\ell}(s),\hat y\rangle \,
\RRR\di{s}.\EEE$$   By \eqref{eq:new-est-strain}, there exists $e\in L^\infty(0,T;L^p(\Omega;\R^{3\times 3})$ such that, for a not relabeled subsequence,
$$
\nabla \overline{y}_\tau \overline{G}_\tau^{-1}\wk^{*} e\quad 
\text{ in } L^{\infty}(0,T;L^p(\Omega;\R^{3\times 3})).
$$
On the other hand, the convergences \eqref{eq:conv3} and
\eqref{eq:conv1bis} yield \RRR that \EEE $e=\nabla y G^{-1}$. We can hence use
convergence \eqref{eq:FH1} \RRR and the polyconvexity from {\bf(H1)}
in order to \EEE pass to the $\liminf$ in
\eqref{eq:dispo} and obtain  that 
\begin{align}
&\RRR \int_t^{t+\delta } \EEE \!\! \int_{\Omega}W(\nabla
                 y(t)\,G^{-1}(s))\, \UUU \det G(s)\EEE\,\di{x}\, \RRR
                 \di{s} \EEE-\RRR \int_t^{t+\delta } \EEE
                 \!\!\langle{\ell}(s),{y}(s)\rangle\, \RRR\di{s} \EEE
                 \leq \liminf_{\tau\to 0}
  \RRR \int_t^{t+\delta } \EEE \!\!\int_{\Omega}W(\nabla \hat y\,\overline{G}_\tau^{-1}(s))\, \UUU \det \overline{G}_\tau(s)\EEE\,\di{x}\, \RRR\di{s} \EEE-\RRR \int_t^{t+\delta } \EEE \!\!\langle
  {\ell}(s),\hat y\rangle\, \RRR\di{s} \EEE
  \nonumber\\
&\quad \stackrel{\eqref{eq:FH1}}{\leq}\RRR \int_t^{t+\delta } \EEE \!\!\int_{\Omega}W(\nabla \hat y\,{G}^{-1}(s))\, \UUU \det G(s)\EEE\,\di{x}\, \RRR\di{s} \EEE-\RRR \int_t^{t+\delta } \EEE \!\!\langle
  {\ell}(s),\hat y\rangle\, \RRR\di{s} \EEE. \label{passit}
\end{align}\RRR

By applying again Lemma \ref{lemma:ball}, this time for the
choice $F = \nabla \hat y\, G^{-1}(t)$ and $H = G(t)G^{-1}(s)$ and using
the time regularity of $G$ and $\ell$ one proves
that 
\begin{equation}
  \label{eq:ball2}
  \lim_{\delta \to 0} \left(\frac{1}{\delta}\RRR \int_t^{t+\delta }  \!\!\int_{\Omega}W(\nabla \hat y\,{G}^{-1}(s))\,  \det G(s)\,\di{x}\, \di{s} - \frac{1}{\delta}\int_t^{t+\delta }  \!\!\langle
  {\ell}(s),\hat y\rangle\, \di{s} \right) =
                \int_{\Omega}W(\nabla \hat y\,{G}^{-1}(t))\,  \det
                G(t)\,\di{x}\,   -  \langle
  {\ell}(t),\hat y\rangle.
\end{equation}
Estimate \eqref{passit} implies that the function  
\begin{equation}
t \in (0,T)\mapsto \int_{\Omega}W(\nabla
                 y(t)\,G^{-1}(t))\,  \det G(t)\,\di{x} - \langle{\ell}(t),{y}(t)\rangle\label{eq:function}
\end{equation} \RRR
is integrable. Choose now $t\in(0,T)$ to be one of its Lebesgue
points. By using \eqref{eq:ball2} one can pass to
the limit as $\delta \to 0$ in \eqref{passit} and deduce that 
$$\int_{\Omega}W(\nabla y(t)\,{G}^{-1}(t))\,  \det
                G(t)\,\di{x}\,   -  \langle
  {\ell}(t),y(t)\rangle \leq \int_{\Omega}W(\nabla \hat y\, {G}^{-1}(t))\,  \det
                G(t)\,\di{x}\,   -  \langle
  {\ell}(t),\hat y\rangle.$$\RRR
As $\hat y\in \mathcal{Y}$ is arbitrary, minimality  
\eqref{eq:problemP} follows. \EEE

\section{Proof of Theorem 2: Optimal
  control}\label{sec:control}

We now turn to the existence proof of optimal controls $\mu^*$ and
optimal pairs $(y^*,G^*) \in S(\mu^*)$ solving problem
\eqref{eq:oc}.

The first step is to check that the solution operator
$S$ is well-defined. This amounts in proving the existence of
nutrient-driven morphoelastic solutions  for given $\mu \in
L^p(0,T;W^{1,p}(\Omega))$, and can be ascertained by extending the argument
of Theorem \ref{thm:main}. Indeed, by resorting again to a time
discretization with constant time step $\tau$, starting from $(y_0,G_0)=(y^0,G^0) \in \mathcal{Y} \times 
\mathcal{G}_p$, one solves for $\{ (y_i,G_i) \}_{i=1}^N \in \mathcal{Y}^N \times 
\mathcal{G}_p^N$ such that, for $i=1,\dotsc,N$, 
\begin{align}
	&\label{eq:discrete-cgrowth1c} y_i\in \text{argmin}_{y \in \mathcal{Y}}\,\left\{\int_{\Omega} W(\nabla y \, G_{i}^{-1})\, \UUU \det G_i\EEE\,\di x- \langle \ell_i,y\rangle
	\right\}, \\
	&\label{eq:discrete-cgrowth2c} G_i = \exp \big(\tau \,  M(G_{i-1}, (K_\tau {\nabla y})_{i-1}, \mu_i) \big) \, G_{i-1} \quad \text{ a.e. in } \Omega, 
\end{align}
where now the additional datum $\mu_i$ is defined as  
\begin{align*}
	\mu_i&:=\frac{1}{\tau} \int_{(i-1)\tau}^{i\tau} \mu(t) \, \di t \in W^{1,p}(\Omega) \quad \text{ for } i=1,\dotsc,N.
\end{align*}
Note here that $\| \overline \mu_\tau \|_{L^p(0,T;W^{1,p}(\Omega))}
\leq \|  \mu \|_{L^p(0,T;W^{1,p}(\Omega))}$ and $\overline \mu_\tau
\to \mu$ in $ L^p(0,T;W^{1,p}(\Omega))$ as $\tau \to 0$.
The existence of a solution to
\eqref{eq:discrete-cgrowth1c}-\eqref{eq:discrete-cgrowth2c} can be obtained by simply
adapting the argument of Lemma  \ref{lemma:exist-y}. The only
modification is required in the estimate on $\nabla G_i$, which now hinges
on  {\bf (H7)} and reads
 \begin{align*}
	\|\nabla G_i \|_{L^p(\Omega; \R^{3\times 3 \times 3})} 
	& \leq  \| G_{i-1}\|_{L^\infty} \| \nabla \exp \big( \tau \,  M(G_{i-1}, (K_\tau {\nabla y})_{i-1}, \mu_i)\big)\|_{L^p(\Omega; \R^{3\times 3 \times 3})}   \\ 
	& \quad +  \| \exp\big(\tau \,  M(G_{i-1}, (K_\tau {\nabla y})_{i-1}, \mu_i) \big) \|_{L^\infty}  
	\| \nabla G_{i-1} \|_{L^p(\Omega; \R^{3\times 3 \times 3})}  \\[0.2cm]
	& \leq \| G_{i-1}\|_{L^\infty} \| \nabla \exp \big( \tau \,  M(G_{i-1}, (K_\tau {\nabla y})_{i-1}, \mu_i)\big)\|_{L^p(\Omega; \R^{3\times 3 \times 3})} \\
	& \quad + C \|  G_{i-1} \|_{W^{1,p}(\Omega; \R^{3\times 3})}. 
 \end{align*}
Here and in the following, we use the symbol $C$ to indicate a positive
 constant, possibly depending on data but not on $\mu$ nor on $\tau$. 
 The actual value of $C$ can change from line to line.

By arguing as in \eqref{eq:est-grad-exp}, we get 
\begin{equation*}
|\nabla  \exp \big(\tau \,  M(G_{i-1}, (K_\tau {\nabla y})_{i-1},\mu_i)\big) | 
 \leq 3\tau \exp \big(  \tau \Vert M\Vert_{L^{\infty}} \big) 
	\Vert M\Vert_{W^{1,\infty}} 
	\Big( |\nabla G_{i-1} | + |\nabla (K_\tau {\nabla y})_{i-1}| + |\nabla \mu_i|  \Big),
\end{equation*} 
so that, using once more {\bf (H7)} we have 
\begin{align*}
&\|  \nabla  \exp \big(\tau \,  M(G_{i-1}, (K_\tau {\nabla y})_{i-1}, \mu_i)\big)  \|_{L^p(\Omega; \R^{3\times 3 \times 3})} 
\\
&\quad \leq C \tau (\|  G_{i-1} \|_{W^{1,p}(\Omega; \R^{3\times 3})}+  \|\nabla (K_\tau {\nabla
	y})_{i-1}\|_{L^p(\Omega; \R^{3\times 3 \times 3})}  + \| \mu_i \|_{W^{1,p} (\Omega)}).
    \end{align*} 
Owing to  \eqref{eq:discrete-conv} and {\bf (H4)} along with the fact that 
$y_{j}\in \mathcal{Y}$ for $j=0,\dots,i-1$, 
\begin{equation}
\label{eq:est-grad-conv-c}
	\|\nabla (K_\tau {\nabla y})_{i-1}\|_{L^p(\Omega; \R^{3\times 3\times 3})}\leq \sum_{j=0}^{i-1}\tau  | \kappa_j |  \|{\nabla y}_{i-1-j}\|_{L^p(\Omega; \R^{3\times 3})}\|\nabla \phi\|_{L^1(\R^3; \R^{3})}\leq C,
\end{equation}
which, together with the facts that $\overline \mu_\tau$ is bounded in 
$L^p(0,T;W^{1,p}(\Omega))$ independently of $\tau$ 
and that $G_{i-1} \in \mathcal{G}_p$, implies that  $G_i\in W^{1,p}(\Omega; \R^{3\times
  3})$.

In view of passing to the limit in the time discretization, a priori
estimates independent of $\tau$ have to be provided. 
The extra $\mu$-dependence of the growth-rate function $M$ has no
influence on estimates \eqref{eq:hatG-winfty}-\eqref{eq:unif-y} and
\eqref{eq:bd-disc-conv1}-\eqref{eq:bd-disc-conv2}, which can be
readily obtained as in Proposition \ref{prop:apriori}. As regards the
estimate on $\overline G_\tau$, one is asked to deal with an extra term
featuring $\nabla \mu_i  $. In particular, subtracting
$G_{i-1}$  from both sides of \eqref{eq:discrete-cgrowth2c}, dividing by $\tau$, and taking the gradient we find 
\begin{equation}
	\label{eq:ing-c1}
	\nabla \left( \frac{G_i - G_{i-1}}{\tau} \right) = \nabla(
        E_{i-1}  G_{i-1})  = ( G_{i-1}^\top \nabla E_{i-1}^\top )^t+ E_{i-1}\nabla G_{i-1} 
\end{equation} 
a.e. in $\Omega$, for $i=1,\dotsc,N$, 
where 
$$
	E_{i-1} := \frac{ \exp \big( \tau \, M(G_{i-1}, (K_\tau
          {\nabla y})_{i-1},\mu_i)\big) - {\rm Id}}{\tau}. 
$$
We can hence control the $L^p$ norm as follows
\begin{equation}
\label{eq:ing-c20}
\left \Vert 	\nabla \left( \frac{G_i - G_{i-1}}{\tau} \right) \right \Vert_{L^p(\Omega; \R^{3\times 3 \times 3})} 
\leq \Vert G_{i-1}^\top \nabla E_{i-1}^\top \Vert_{L^p(\Omega; \R^{3\times 3 \times 3})}  + 
\Vert E_{i-1} \nabla G_{i-1} \Vert_{L^p(\Omega; \R^{3\times 3 \times 3})} .
\end{equation}
From the properties of the matrix exponential, arguing as in
\eqref{eq:est-grad-exp}, we get that 
\begin{equation*}
|\nabla E_{i-1}|  \leq 3 \exp \big(  \tau \|M\|_{L^{\infty}} \big)  \Vert M \Vert_{W^{1,\infty}} 
\Big( |\nabla G_{i-1} | + |\nabla  (K_\tau {\nabla y})_{i-1}| + | \nabla \mu_i |  \Big)
\end{equation*}
and hence that 
\begin{align}
&\Vert \nabla E_{i-1}  \Vert_{L^p(\Omega; \R^{3\times 3 \times 3})}  \nonumber \\
&\quad \leq C\exp \big(  \tau \|M\|_{L^{\infty}} \big)  \Vert M \Vert_{W^{1,\infty}} \nonumber\\
& \qquad  \times	\Big(   \Vert \nabla G_{i-1}  \Vert_{L^p(\Omega; \R^{3\times 3 \times 3})}   
+  \Vert \nabla  (K_\tau {\nabla y})_{i-1}  \Vert_{L^p(\Omega; \R^{3\times 3 \times 3})}  
+  \Vert  \nabla \mu_i  \Vert_{L^p(\Omega; \R^3)}  \Big) \nonumber \\[0.2cm]
& \quad \leq C\exp \big(  \tau \|M\|_{L^{\infty}} \big)  \Vert M \Vert_{W^{1,\infty}} \nonumber \\
& \qquad  \times	\Big(   \Vert \nabla G_{i-1}  \Vert_{L^p(\Omega; \R^{3\times 3 \times 3})}   
+  \sum_{j=0}^{i-1}\tau  | \kappa_j |  \|{\nabla y}_{i-1-j}\|_{L^p(\Omega; \R^{3\times 3})}\|\nabla \phi\|_{L^1(\R^3; \R^{3})}
+  \Vert \mu_i  \Vert_{W^{1,p}(\Omega)}   \Big)  \nonumber \\[0.2cm]
&\quad \leq  C \left( 1+
\Vert \nabla G_{i-1} \Vert_{L^p(\Omega; \R^{3\times 3 \times 3})} + \Vert \mu_i  \Vert_{W^{1,p}(\Omega)}   \right) 
\label{eq:ing-c30}, 
\end{align}
for all $i =1,\dotsc,N$, where we have also used that $\| y_j
\|_{W^{1,p}(\Omega;\R^3)}$ is bounded independently of $\tau$, as well as that one
can control
\begin{equation}
  \label{eq:est-grad-conv-c2}
	  \sum_{j=0}^{i-1}\tau  | \kappa_j |  \|{\nabla y}_{i-1-j}\|_{L^p(\Omega; \R^{3\times 3})}\|\nabla \phi\|_{L^1(\R^3; \R^{3})}\leq C .
\end{equation} 
On the other hand, similar calculations as in 
\eqref{eq:ing0} shows that 
$$
\Vert  E_{i-1}  \Vert_{L^\infty}  \leq \frac{1}{\tau}   \left(\exp(\tau \|M\|_{L^\infty})-1\right)
\leq C. 
$$
Going back to \eqref{eq:ing-c20} and using the H\"older Inequality, \eqref{eq:ing-c30}, and the fact that 
$\| G_{i-1} \|_{L^\infty}$ is bounded independently of $\tau$, we deduce that 
\begin{align*}
&\left \Vert 	\nabla \left( \frac{G_i - G_{i-1}}{\tau} \right) \right \Vert_{L^p(\Omega; \R^{3\times 3 \times 3})} 
\\
& \quad \leq \Vert G_{i-1} \Vert_{L^\infty}  \Vert  \nabla E_{i-1}  \Vert_{L^p(\Omega; \R^{3\times 3 \times 3})}  
+ \Vert E_{i-1} \Vert_{L^\infty} 
\Vert  \nabla G_{i-1} \Vert_{L^p(\Omega; \R^{3\times 3 \times 3})} \\[2mm]
& \quad \leq C \left( 1+
\Vert \nabla G_{i-1} \Vert_{L^p(\Omega; \R^{3\times 3 \times 3})} + \Vert \mu_i  \Vert_{W^{1,p}(\Omega)}   \right)  \\
& \quad \leq C \left( 1+ \Vert \nabla G_{0} \Vert_{L^p(\Omega; \R^{3\times 3 \times 3})}+
 \sum_{j=1}^{i-1} \tau \left \Vert 	\nabla \left( \frac{G_j - G_{j-1}}{\tau} \right) \right \Vert_{L^p(\Omega; \R^{3\times 3 \times 3})} 
 + \Vert \mu_i  \Vert_{W^{1,p}(\Omega)}   \right) .  
\end{align*}
By taking the $p$-power, applying the Discrete Gronwall Lemma, and
recalling that $\overline \mu_\tau$ is bounded in 
$L^p(0,T;W^{1,p}(\Omega))$ independently of $\tau$ we conclude that
 $\nabla \widehat G_\tau'$ is bounded in $L^p(0,T;L^p(\Omega;
\R^{3\times 3 \times 3}))$, independently of $\tau$. Consequently, 
it follows that $\widehat G_\tau$ is bounded in $W^{1,p}(0,T;W^{1,p}(\Omega;
\R^{3\times 3}))$ and that $\overline G_\tau$ is bounded in $L^\infty(0,T;W^{1,p}(\Omega;
\R^{3\times 3}))$, both independently of $\tau$. 

One can now extract   not relabeled subsequences and pass to the limit
as $\tau \to 0$, following the very argument of Subsection
\ref{subs:limit}. Note nonetheless that the convergence of $\widehat
G_\tau$ and $\overline G_\tau$ is slightly weaker, namely,
\begin{align} 
\label{eq:conv1-c}\widehat{G}_\tau\wstarto G\quad
&\text{ in } W^{1,\infty}(0,T;L^\infty(\Omega;\R^{3\times 3}))\cap
                 W^{1,p}(0,T;W^{1,p}(\Omega;\R^{3\times 3})),\\
	\label{eq:conv2-c}
	\widehat{G}_\tau\to G\quad
	&\text{ in } C([0,T];C(\overline{\Omega};\R^{3\times 3})) ,\\
\label{eq:conv3-c} \overline{G}_\tau\to G\quad
	&\text{ in } L^\infty ((0,T)  \times \Omega ;\R^{3\times 3}) \cap L^p(0,T;L^p(\Omega; \R^{3 \times 3})), \\
\label{eq:conv-c}\overline{G}_\tau\wstarto G\quad 
	&\text{ in } L^\infty(0,T;W^{1,p}(\Omega;\R^{3\times 3})),
\end{align} 
where we have also used the Aubin-Lions Lemma.  As $\overline \mu_\tau \to
\mu$ in $L^p(0,T;W^{1,p}(\Omega))$ and  the growth-rate function
$M$ is Lipschitz continuous with respect to $\mu$  from {\bf
  (H7)}, we readily check again that the limit $(y,G)$ of
time-discrete solutions is a nutrient-driven solution in the sense of
Definition \ref{nutrient-driven}.

The above argument shows that the solution operator
\begin{align*}
&S: L^p(0,T;W^{1,p}(\Omega)) \\
&\quad \to L^\infty(0,T;W^{1,p}(\Omega; \R^3)) 
 \times W^{1,\infty}(0,T;L^\infty(\Omega; \R^{3\times 3}))
\cap W^{1,p}(0,T;W^{1,p}(\Omega; \R^{3\times 3}))
\end{align*}
 representing the
set $S(\mu)$ of all nutrient-driven solutions $(y,G)$ given $\mu$ is
well-defined and bounded. 

In order to prove Theorem \ref{thm:control}, let now $\mu_k\in
\mathcal A$ and $(y_k,G_k)\in S(\mu_k)$ with
$$
J(y_k,G_k,\mu_k) \to \inf_{\mu\in \mathcal A}\big\{J(y,G,\mu) \ : \
  (y,G)\in S(\mu) \big\} \geq 0.
$$
Since $\mathcal A$ is bounded in $L^p(0,T;W^{1,p}(\Omega))$ and compact
in $L^1((0,T)\times \Omega)$ by {\bf (H8)}, one can find $\mu^*\in
\mathcal A$ and pass to a not
relabeled subsequence such that $\mu_k \to \mu^*$ a.e. Moreover, the
boundedness of $S$ implies that $(y_k,G_k)$ are uniformly bounded in 
$L^\infty(0,T;W^{1,p}(\Omega; \R^3))  \times W^{1,\infty}(0,T;L^\infty(\Omega; \R^{3\times 3}))
\cap W^{1,p}(0,T;W^{1,p}(\Omega; \R^{3\times 3}))$. Hence, by
extracting again (without relabeling) we get that
\begin{align*}
  y_k \wstarto y^*\quad 
  &\text{ in } L^\infty(0,T;W^{1,p}(\Omega;\R^{3\times 3})),\\
 G_k \wstarto G^*\quad 
  &\text{ in } W^{1,\infty}(0,T;L^\infty(\Omega; \R^{3\times 3})) \cap W^{1,p}(0,T;W^{1,p}(\Omega; \R^{3\times 3})).
\end{align*}
The latter implies in particular that $G_k \to G^*$ in $C([0,T] \times \overline
\Omega ;\R^{3\times 3})$. 
Moreover, since $\det G_k$ is a.e. bounded below by a positive constant
independently of $k$ one has that  
$$
(G_k)^{-1} \to (G^*)^{-1} \quad \text{in } C([0,T] \times \overline
\Omega;\R^{3\times 3})),
$$
as well. These convergences are enough to pass to the limit in
relations \eqref{eq:problemPn}-\eqref{eq:initialn} and obtain that
the limiting $(y^*,G^*)$ belongs to $S(\mu^*)$. \RRR In particular,
the almost-everywhere-in-time minimality \eqref{eq:problemPn} follows
along the same lines as in the proof of Theorem~1, see Subsection
\ref{subs:limit}. \EEE

On the other hand, due to the lower semicontinuity of $J$ from {\bf
  (H9)} we get that
\begin{align*}
  J(y^*,G^*,\mu^*) \leq \liminf_{k\to \infty} J(y_k,G_k,\mu_k) = \min_{\mu\in \mathcal A}\big\{J(y,G,\mu) \ : \
  (y,G)\in S(\mu) \big\}.
\end{align*}
In particular, $\mu^*$ is an optimal control and $(y^*,G^*)\in S(\mu^*)$
 is the corresponding optimal state.

\section{Proof of Theorem 3: Nutrient-morphoelastic existence}
\label{sec:chemotaxis}
 
In this section, we prove the existence of nutrient-morphoelastic
solutions, namely, trajectories $t\in [0,T]  \mapsto (y(t),G(t),\mu(t))\in \mathcal{Y} 
\times \mathcal{G}_p\times \mathcal{M}$ satisfying \eqref{eq:problemcP}-\eqref{eq:initialc-1}. Using again a time discretization with 
constant time step $\tau$, starting from  $(y_0,G_0,\mu_0)=(y^0,G^0,\mu^0) \in \mathcal{Y} \times 
\mathcal{G}_p\times \mathcal{M}$, we look  for $\{ (y_i,G_i,\mu_i) \}_{i=1}^N \in \mathcal{Y}^N \times 
\mathcal{G}_p^N\times \mathcal{M}^N$ fulfilling 
\begin{align}
	&\label{eq:discrete-growth1cm} 
	y_i\in \text{argmin}_{y \in \mathcal{Y}}\,\left\{\int_{\Omega} W(\nabla y \, G_{i}^{-1})\, \UUU \det G_i\EEE\,\di x- \langle \ell_i,y\rangle
	\right\}, \\
	&\label{eq:discrete-growth2cm} 
	G_i = \exp \big(\tau \,  M(G_{i-1}, (K_\tau {\nabla y})_{i-1}, \mu_{i-1}) \big) \, G_{i-1} \quad \text{ a.e. in } \Omega, \\
	&\label{eq:discrete-growth3cm} 
	\frac{\mu_i-\mu_{i-1}}{\tau}- \nu \Delta \mu_i = h_i - H( (\kappa \ast_{\tau} y)_{i-1}) \quad \text{ a.e. in } \Omega, \\
	& \label{eq:discrete-growth4cm} 
	\mu_i= \mu_{{\rm D},i}  \quad \text{ a.e. on } \partial \Omega
\end{align}
for $i =1,\dotsc,N$. In \eqref{eq:discrete-growth3cm}-\eqref{eq:discrete-growth4cm}, we have set
$$
h_i:= \frac{1}{\tau} \int_{(i-1)\tau}^{i\tau} h(s) \,\di{s} \in L^p(\Omega) \quad \text{ for } i =1,\dotsc,N, 
$$
and 
$$
\mu_{{\rm D},i}:= \frac{1}{\tau} \int_{(i-1)\tau}^{i\tau} \mu_{\rm D}(s) \,\di{s} \in W^{2,p}(\Omega) \quad \text{ for } i =1,\dotsc,N, 
$$
and we recall that 
$$
(\kappa \ast_{\tau} y)_{i-1}(x):= \sum_{j=0}^{i-1} \tau \,  \kappa_j \, y_{i-1-j}(x) \quad \text{ for a.e. } x \in \Omega, \text{ for }   i =1,\dotsc,N, 
$$
where $\kappa_i:= \kappa(t_i)$,  $i =0,\dotsc,N$. 

Let us first prove that the scheme \eqref{eq:discrete-growth1cm}-\eqref{eq:discrete-growth4cm} admits a solution. 
Using the same arguments as in Lemma \ref{lemma:exist-y}, it follows that for $i =0,\dotsc,N-1$, given $G_i \in \mathcal{G}_p$, 
there exists $y_i \in \mathcal{Y}$ solving \eqref{eq:discrete-growth1cm}. By quite similar arguments to those in the 
beginning of Section \ref{sec:control}, it is readily verified that for $i =1,\dotsc,N$, given $(y_{i-1}, G_{i-1}, \mu_{i-1}) \in  
\mathcal{Y} \times \mathcal{G}_p\times \mathcal{M}$, there exists a solution $G_i\in \mathcal{G}_p$ to \eqref{eq:discrete-growth2cm}. 
It remains to check that for $i =1,\dotsc,N$, given $(y_{i-1}, \mu_{i-1}) \in  \mathcal{Y} \times \mathcal{M}$, there exists a solution 
$\mu_i \in \mathcal{M}$ to  \eqref{eq:discrete-growth3cm}-\eqref{eq:discrete-growth4cm}. Letting 
$$
F_i := \tau h_i -\tau H ((\kappa \ast_{\tau} y)_{i-1})  + \mu_{i-1} \in L^p(\Omega),
$$
an application of \cite[Theorem 2.4.2.5]{grisvard} shows that the problem  
\begin{equation}
	\label{eq:discrete-growth5cm} 
\mu_i - \nu \tau \Delta \mu_i = F_i \: \text{ a.e. in } \Omega, \quad 	\mu_i= \mu_{{\rm D},i}  \:  \text{ a.e. on } \partial \Omega, 
\end{equation}
has a unique solution $\mu_i \in \mathcal{M}$. 

We next perform some a-priori estimates. The additional dependence of the growth-rate function $M$ on $\mu$ has no impact 
on the estimates \eqref{eq:hatG-winfty}-\eqref{eq:unif-y} and \eqref{eq:bd-disc-conv1}-\eqref{eq:bd-disc-conv2}. These 
estimates can be obtained as in Proposition \ref{prop:apriori}. For the estimates on $\bar{\mu}_{\tau}$ and $\hat{\mu}_{\tau}$, we 
deduce by applying \cite[Theorems 2.3.3.6 \& 1.5.1.2]{grisvard} and using \eqref{eq:discrete-growth5cm} that 
\begin{align}
\label{eq:apriori-cm} 
	\Vert \mu_i \Vert_{W^{2,p}(\Omega)} &\leq C \left( \Vert \mu_i - \nu \tau \Delta \mu_i \Vert_{L^{p}(\Omega)} 
	+ \Vert \mu_i \vert_{\partial \Omega} \Vert_{W^{2-1/p,p}(\partial \Omega)} \right) \nonumber \\
	&\leq C \left( \Vert F_i \Vert_{L^{p}(\Omega)} + \Vert \mu_{{\rm D},i}\Vert_{W^{2,p}(\Omega)} 
	\right) 
\end{align}
for every $i =1,\dotsc,N$. Here and in the following, the symbol $C$ stands for a positive constant, possibly depending on the data 
but not on $\tau$ and varying from line to line. By letting $\tilde{\mu}_i:= \mu_i -\mu_{{\rm D},i}$ and 
$$
\tilde{F}_i:= h_i - H((\kappa \ast_\tau y)_{i-1}) - \frac{\mu_{{\rm D},i}- \mu_{{\rm D},i-1}}{\tau}  + \nu \Delta \mu_{{\rm D},i} \in L^p(\Omega), 
$$
equations \eqref{eq:discrete-growth3cm} and \eqref{eq:discrete-growth4cm} read
\begin{align}
&\label{eq:discrete-growth6cm} 
\frac{\tilde{\mu}_i - \tilde{\mu}_{i-1}}{\tau}- \nu \Delta \tilde{\mu}_i = \tilde{F}_i   \quad \text{ a.e. in } \Omega, \\
& \label{eq:discrete-growth7cm} 
\tilde{\mu}_i= 0 \quad \text{ a.e. on } \partial \Omega,
\end{align}
for $i =1,\dotsc,N$. We multiply \eqref{eq:discrete-growth6cm} by the function $\tau |\tilde{\mu}_i|^{p-2} \tilde{\mu}_i$ and integrate over $\Omega$ 
to obtain that 
\begin{align}
\label{eq:apriori-cm1} 
\int_{\Omega}  \left(\tilde{\mu}_i - \tilde{\mu}_{i-1} \right) |\tilde{\mu}_i|^{p-2} \tilde{\mu}_i \,\di{x}  - \nu \tau \int_{\Omega}  \left(\Delta \tilde{\mu}_i \right) |\tilde{\mu}_i|^{p-2} \tilde{\mu}_i \,\di{x} 
= \tau \int_{\Omega}  \tilde{F}_i  \,  |\tilde{\mu}_i|^{p-2} \tilde{\mu}_i \,\di{x}.
\end{align}
From the convexity of the map $\mu \mapsto |\mu|^p/p$ we obtain
\begin{align}
\label{eq:apriori-cm2} 
\int_{\Omega}  \left(\tilde{\mu}_i - \tilde{\mu}_{i-1} \right) |\tilde{\mu}_i|^{p-2} \tilde{\mu}_i \,\di{x} 
& \geq \frac{1}{p}  \int_{\Omega} |\tilde{\mu}_i|^{p} \,\di{x} - \frac{1}{p}  \int_{\Omega} |\tilde{\mu}_{i-1}|^{p} \,\di{x}.
\end{align}
On the other hand, by applying the H\"older and the Young Inequalities we have  
\begin{align}
	\label{eq:apriori-cm3} 
\tau \int_{\Omega}  \tilde{F}_i  \,  |\tilde{\mu}_i|^{p-2} \tilde{\mu}_i \,\di{x} 
&\leq \tau \left(   \int_{\Omega} |\tilde{F}_i|^{p} \,\di{x} \right)^{1/p} \left(   \int_{\Omega} |\tilde{\mu}_i|^{p} \,\di{x} \right)^{1/q} \nonumber \\
& \leq \tau \left(  \frac{1}{p}  \int_{\Omega} |\tilde{F}_i|^{p} \,\di{x} +  \frac{1}{q}  \int_{\Omega} |\tilde{\mu}_i|^{p} \,\di{x} \right). 
\end{align}
Furthermore, by using the Green Formula and taking into account the boundary condition \eqref{eq:discrete-growth7cm}, we get 
\begin{align}
	\label{eq:apriori-cm4} 
- \nu \tau \int_{\Omega}  &( \Delta \tilde{\mu}_i )  \, |\tilde{\mu}_i|^{p-2} \tilde{\mu}_i \,\di{x} 
= \nu \tau (p-1) \int_{\Omega}  |\tilde{\mu}_i|^{p-2} | \nabla \tilde{\mu}_i|^2 \,\di{x} . 
\end{align}
Thus, by plugging \eqref{eq:apriori-cm2}-\eqref{eq:apriori-cm4} in \eqref{eq:apriori-cm1}, we arrive at 
$$
\frac{1}{p}  \int_{\Omega} |\tilde{\mu}_i|^{p} \,\di{x} - \frac{1}{p}  \int_{\Omega} |\tilde{\mu}_{i-1}|^{p} \,\di{x} 
+  \nu \tau (p-1) \int_{\Omega}  |\tilde{\mu}_i|^{p-2} | \nabla \tilde{\mu}_i|^2 \,\di{x}  
\leq  \frac{\tau}{p}  \int_{\Omega} |\tilde{F}_i|^{p} \,\di{x} +  \frac{\tau}{q}  \int_{\Omega} |\tilde{\mu}_i|^{p} \,\di{x}. 
$$
Summing this inequality from $i=1$ to $m\leq N$ we infer that 
$$
\left( \frac{1}{p} - \frac{\tau}{q} \right) \Vert \tilde{\mu}_m \Vert^p_{L^p(\Omega)} \leq  \frac{1}{p}  \Vert \tilde{\mu}_0 \Vert^p_{L^p(\Omega)} 
+ \frac{\tau}{p} \sum_{i=1}^m \Vert \tilde{F}_i \Vert^p_{L^p(\Omega)} + \frac{\tau}{q} \sum_{i=1}^{m-1} \Vert \tilde{\mu}_i \Vert^p_{L^p(\Omega)}. 
$$
Upon choosing $\tau$ small enough, by applying the Discrete
Gronwall Lemma and owing to assumptions {\bf(H10)} and {\bf(H11)} we deduce that $
\Vert \tilde{\mu}_m \Vert_{L^p(\Omega)} \leq C$ for all
$m=1,\dotsc,N$. Hence, the definition of $\tilde{\mu}_m$, the reverse triangle inequality, and the fact that 
$\overline{\mu_{\rm D}}\,\hspace{-0.05cm}_{,\tau}$ is bounded in $L^\infty(0,T;W^{2,p}(\Omega))$ independently of $\tau$, result in $
\Vert \mu_m \Vert_{L^p(\Omega)} \leq C$ for all $m=1,\dotsc,N$. 
Now, inserting this estimate into \eqref{eq:apriori-cm} and using again {\bf(H10)} and {\bf(H11)}, we find that, for sufficiently small $\tau$, 
$$
\Vert \mu_i \Vert_{W^{2,p}(\Omega)} \leq C \quad \forall i=1,\dotsc,N.
$$
This proves that $\bar{\mu}_{\tau}$ and $\hat{\mu}_{\tau}$ are bounded in $L^\infty(0,T; W^{2,p} (\Omega))$ independently of $\tau$, and, together with 
\eqref{eq:discrete-growth3cm}, {\bf(H10)}, and {\bf(H11)}, it also proves that $\hat{\mu}_{\tau}$ is bounded in $W^{1,\infty}(0,T; L^{p} (\Omega))$ independently of~$\tau$. 

Using the fact that $\bar{\mu}_{\tau}$ is bounded in $L^\infty(0,T;
W^{2,p} (\Omega))$ independently of $\tau$, the same arguments
from Section 
\ref{sec:control}  entail that $\widehat G_\tau$ is bounded in $W^{1,p}(0,T;W^{1,p}(\Omega;
\R^{3\times 3}))$ and that $\overline G_\tau$ is bounded in $L^\infty(0,T;W^{1,p}(\Omega;
\R^{3\times 3}))$, both independently of $\tau$. 

We proceed to show some a-priori bounds for $\overline{(\kappa \ast_\tau  y)}_\tau$ and
$\widehat{(\kappa \ast_\tau  y)}_\tau$. By \eqref{eq:unif-y} and {\bf (H4)}, we have for every $i=1,\dots,N$ 
\begin{align*}
	 \| (\kappa \ast_\tau  y)_i\|_{L^p(\Omega; \R^3)} &\leq \sum_{j=0}^N \tau |\kappa_j|  \| y_{i-j}\|_{L^p(\Omega; \R^{3})}  \\
	&
	\leq \tau (N+1) \|\kappa\|_{L^\infty(0,T)} \|  \overline{y}_\tau\|_{L^{\infty}(0,T;L^p(\Omega;\R^{3}))} \leq C
\end{align*} 
and analogously 
$$
\| \nabla (\kappa \ast_\tau  y)_i\|_{L^p(\Omega; \R^{3 \times 3})} 
\leq \tau (N+1) \|\kappa\|_{L^\infty(0,T)} \| \nabla  \overline{y}_\tau\|_{L^{\infty}(0,T;L^p(\Omega;\R^{3 \times 3}))} \leq C, 
$$
which imply that both $\overline{(\kappa \ast_\tau  y)}_\tau$ and $\widehat{(\kappa \ast_\tau  y)}_\tau$ are bounded in $L^{\infty}(0,T;W^{1,p}(\Omega;\R^{3}))$ independently of 
$\tau$. Moreover for a.e. $x \in \Omega$, 
\begin{align*}
\left|	\frac{(\kappa \ast_\tau y)_i(x)-(\kappa \ast_\tau y)_{i-1}(x)}{\tau} \right| &= \left|  \kappa_0 y_i(x) + \sum_{j=1}^i (\kappa_j-\kappa_{j-1})  y_{i-j}(x)  \right| \\
&\leq \|\kappa\|_{L^\infty(0,T)}  |y_i(x) | + \|\kappa'\|_{L^1(0,T)} | y_{i-j}(x)|, 
\end{align*}
so that 
$$
\left \Vert \frac{(\kappa \ast_\tau y)_i -(\kappa \ast_\tau y)_{i-1}}{\tau} \right \Vert_{L^p(\Omega; \R^3)} \leq C \left( \|\kappa\|_{L^\infty(0,T)} + 
 \|\kappa'\|_{L^1(0,T)} \right) \|  \overline{y}_\tau\|_{L^\infty(0,T;L^p(\Omega;\R^{3}))}. 
$$
This, along with \eqref{eq:unif-y} and  {\bf (H4)}, yields that $\widehat{(\kappa \ast_\tau  y)}_\tau'$ is bounded in $L^{\infty}(0,T;L^{p}(\Omega;\R^{3}))$ independently of 
$\tau$. 

We can now pass to the limit as $\tau \to 0$. Following the arguments in Sections \ref{subs:limit} and \ref{sec:control}, we can extract subsequences, without relabeling, 
in such a way that 
\begin{align*} 
\overline{ y}_\tau	\wstarto y \quad 
&\text{ in } L^{\infty}(0,T;W^{1,p}(\Omega;\R^{3})), \\	
\widehat{G}_\tau \wstarto G\quad
&\text{ in } W^{1,\infty}(0,T;L^\infty(\Omega;\R^{3\times 3}))\cap
W^{1,p}(0,T;W^{1,p}(\Omega;\R^{3\times 3})),\\
 \overline{G}_\tau\to G\quad
&\text{ in } L^\infty ((0,T)  \times \Omega ;\R^{3\times 3}) \cap L^p(0,T;L^p(\Omega; \R^{3 \times 3})), \\
\overline{(K_\tau {\nabla y})}_\tau \to K\nabla y \quad
&\text{ in } L^\infty ((0,T)  \times \Omega ;\R^{3\times 3}) \cap L^p(0,T;L^p(\Omega; \R^{3 \times 3})),
\end{align*} 
and 
\begin{align*} 
\widehat{ \mu}_\tau	\wstarto \mu \quad 
&\text{ in } L^{\infty}(0,T;W^{2,p}(\Omega)) \cap W^{1,\infty}(0,T;L^p(\Omega)), \\	
\overline{\mu}_\tau \to  \mu \quad
&\text{ in } L^{\infty}(0,T;L^p(\Omega)),\\
\overline{\mu}_\tau \wstarto  \mu \quad
&\text{ in } L^\infty (0,T; W^{2,p}(\Omega )), \\
\overline{(\kappa \ast_\tau y)}_\tau  \to \kappa \ast y \quad
&\text{ in } L^{\infty}(0,T;L^p(\Omega;\R^3)).
\end{align*} 
Since the functions $M$ and $H$ are Lipschitz continuous according to  {\bf (H7)} and  {\bf (H11)}, and since $\overline{h}_\tau \to h$ in $L^{\infty}(0,T;L^p(\Omega))$, 
we easily check that the limit $(y,G,\mu)$ is a nutrient-morphoelastic
solution in the sense of Definition \ref{nutrient-morphoelastic}. \RRR
Once again, 
the almost-everywhere-in-time minimality \eqref{eq:problemcP} follows
by adapting the argument of Subsection
\ref{subs:limit}. \EEE

\section*{Acknowledgements}
E.D. acknowledges support from the Austrian Science Fund (FWF) through projects F\,65,  I\,4052, V\,662, and Y\,1292, as well as from BMBWF through the OeAD-WTZ project CZ04/2019. K.N. is partially supported by the Austrian Science Fund (FWF) project  F\,65. U.S. is supported by the Austrian Science Fund (FWF)
projects F\,65, W\,1245, I\,4354, I\,5149, and P\,32788 and by the
OeAD-WTZ project CZ 01/2021. \RRR The Authors are indebted to Matthias
Liero and Willem van Oosterhout for some comments on a previous version of the manuscript. \EEE

\end{document}